%% file: manuscript.tex
\RequirePackage[reqno]{amsmath} 
\documentclass[11pt,reqno]{amsart}

\usepackage[foot]{amsaddr}

\usepackage{mypub}
\usepackage{subfiles}

\usepackage{geometry}
\begin{document}

\title[Global existence for a nonlinear wave equation]{Global
  existence of a  nonlinear wave equation arising from Nordström's
  theory of gravitation}


  \author[U.~Brauer]{Uwe Brauer}

  \address{%
    Uwe Brauer Departamento de Matemática Aplicada\\ Universidad
    Complutense Madrid 28040 Madrid, Spain} \email{oub@mat.ucm.es}

  \thanks{U.~B.~gratefully acknowledges support from Grant PID2019-103860GB-I00
    MINECO, Spain and UCM-GR17-920894.}

\author[L.~Karp]{Lavi Karp}

  \address{%
    Lavi Karp
    Department of Mathematics\\ ORT Braude College\\
    P.O. 
    Box 78, 21982 Karmiel\\ Israel}

  \email{karp@braude.ac.il}

  \subjclass{35Q31 (35B30 35L45 35L60)}

  \keywords{Global existence, Nordström theory, semi-linear wave equation,
homogeneous Sobolev spaces, hyperbolic symmetric systems, energy estimates}

\begin{abstract}

  We show global existence of classical solutions for the nonlinear
  Nordström theory with a source term and a cosmological constant
  under the assumption that the source term is small in an appropriate
  norm, while in some cases no smallness assumption on the initial
  data is required.
  In this theory, the gravitational field is described by a single
  scalar function that satisfies a certain semi-linear wave equation.
  We consider spatial periodic deviation from the background metric,
  that is why we study the semi-linear wave equation on the
  three-dimensional torus $\setT^3$ in the Sobolev spaces $H^m(\setT^3)$.
  We apply two methods to achieve the existence of global solutions,
  the first one is by Fourier series, and in the second one, we write
  the semi-linear wave equation in a non-conventional way as a
  symmetric hyperbolic system.
  We also provide results concerning the asymptotic behavior of these
  solutions and, finally, a blow-up result if the conditions of our
  global existence theorems are not met.
\end{abstract}

\maketitle{}

\tableofcontents

\def\biblio{}

\subfile{section1-intro.tex}

\subfile{section2-preliminaries.tex}

\subfile{section3-fourier.tex}

\subfile{section4-sh.tex}

\subfile{section5-blowup.tex}

\subfile{section6-positivity.tex}
\bibliography{}

\end{document}

%% file: section1-intro.tex
\section{Introduction}
\label{sec:nordstr-euler-syst}

The purpose of the work presented here is to prove global existence
and uniqueness of classical solutions and its asymptotic behavior of a
semi-linear wave equation with damping terms.
This wave equation arises in the context of the nonlinear Nordström
theory of gravity, which we shall describe in what follows.

The first fully relativistic, consistent, theory of gravitation was a
scalar theory developed by Nordström
\cite{nordstrom13:_zur_theor_gravit_stand_relat
}, where the gravitational field is described by a nonlinear
hyperbolic equation for the scalar field $\phi$.
Although the theory is not in agreement with observations  it
provides, due to its nonlinearity, some interesting mathematical
challenges. 
Surprisingly, this theory has never been mathematically investigated,
although its linear version coupled to the Euler equations has been
studied by Speck
\cite{Speck_09
}
and coupled to the Vlasov equation by Calogero
\cite{calogero03:_spher
} and others~\cite{Felix_Antonio_Calogero_2014
}, \cite{Fajman_Jeremie_Jacques-2021
}, \cite{Wang_2021
}, \cite{Calogero_Rein_2003
} and \cite{Calogero_Rein_2004
}.

We follow here the geometric reformulation provided by Einstein-Fokker
\cite{einstein14:_nords_gravit_stand_differ
} and will use the Euler equations as a matter model. 
See also Straumann~\cite[Chap.~2.]{STRA2
} for a modern representation of that theory.
The basic idea of this theory is that the physical metric
$ g_{\alpha\beta}$ is related to the Minkowski metric
$\eta_{\alpha\beta}$ by the following conformal transformation.
\begin{equation}
  \label{eq:Euler-Nordstrom:1}
   g_{\alpha\beta}=\phi^2\eta_{\alpha\beta}, 
\end{equation}
where $\eta_{\alpha\beta}=\text{diag}(-1,1,1,1)$.
The matter is described by an energy-momentum tensor, which in the
case of a perfect fluid takes the form
\begin{equation}
  \label{eq:section1-intro:6}
   T^{\alpha\beta}= \left(  \epsilon + p
  \right) u^{\alpha} u^{\beta}+  p  
g^{\alpha\beta},
\end{equation}
where $\epsilon$ denotes the energy density, $p$ the pressure 
and $u^{\alpha}$ is the unit timelike vector which satisfies
\begin{equation}
  \label{eq:section1-intro:7}
   g_{\alpha\beta} u^{\alpha} u^{\beta}=-1.
\end{equation}
The field equations, as proposed by Einstein and Fokker takes the
following form
\begin{equation}
  \label{eq:Euler-Nordstrom:5}
   R =  T,
\end{equation}
here we set the relevant constants to one, the Ricci scalar is denoted
by $R$ and the trace of the fundamental energy tensor by
$T=g_{\alpha\beta}T^{\alpha\beta}$.
Using equation \eqref{eq:Euler-Nordstrom:1}, the Ricci scalar takes the form
\begin{equation}
\label{eq:section1-intro:2}
  R = -6  \frac{\Box\phi}{\phi^3}, \qquad\Box 
\eqdef\eta^{\alpha\beta}\partial_{\alpha}\partial_{\beta}.
\end{equation}

While the Euler equations take the form
\begin{equation}
  \label{eq:Euler-Nordstrom:6}
   \nabla_{\alpha} T^{\alpha\beta}=0,
\end{equation}
where $ \nabla_{\alpha}$ is the covariant derivative associated
with $g_{\alpha\beta}$.
Combining equations \eqref{eq:Euler-Nordstrom:5}, \eqref{eq:section1-intro:2} and 
\eqref{eq:Euler-Nordstrom:6},  the Euler-Nordström  system takes the
following form
\begin{subequations}
\begin{align}
\label{eq:Euler-Nordstrom:7}
& \Box\phi              =- \frac{1}{6} T\phi^{3} \\
\label{eq:section1-intro:4}  
 &   \nabla_{\alpha} T^{\alpha\beta}  =0.
\end{align}
\end{subequations}
\begin{rem}[Different form of the field equation]
\label{rem:section1-intro:1}
We  want to point out that it is possible to consider a slightly
different conformal transformation (see for example
\cite{calogero03:_spher
} or~\cite{Speck_09
}), namely
\begin{equation}
\label{eq:section1-intro:1}
   g_{\alpha\beta}=e^{2\psi}\eta_{\alpha\beta}, 
\end{equation}
which leads to an equivalent nonlinear wave equation
\begin{equation}
  \label{eq:section1-intro:3}
  \Box\psi + \left( \nabla\psi \right)^2 =-\frac{1}{6} e^{2\psi}T.
\end{equation}
\end{rem}

\subsection{The field equations with  cosmological constant and the background 
solutions}
\label{sec:fields-equations}
In what follows we modify the field equation
\eqref{eq:Euler-Nordstrom:7} by adding a term which corresponds to the
cosmological constant $\Lambda$ in General Relativity in the following way,
\begin{equation}
  \label{eq:Nordstrom:3}
  \Box\phi =- \frac{1}{6} T\phi^{3} -\Lambda \phi.
\end{equation}
This choice is motivated by the properties of explicit solutions which
are homogeneous and isotropic, namely that these properties are very
similar to the ones of Euler-Einstein (see e.g.~\cite[Chap.~V]{Choquet-Bruhat_09
}, \cite[Chap.~10]{Rendall_book
}), and Euler-Poisson 
(\cite{Brauer_Rendall_Reula_94
}), which we will discuss below.

We denote an isotropic and homogeneous vacuum background solution by
$\mathring{\phi}$, and for convenience we set $\varkappa^2=\Lambda>0$.
Homogeneity implies that the function $\mathring \phi$ depends just on
$t$, while the fact that the solution describes vacuum leads to the
conclusion that $T\equiv 0$. 
Therefore equation \eqref{eq:Nordstrom:3} reduces to
\begin{equation}
  \label{eq:section1-intro:1B}
  -\frac{d^2}{dt^2}\mathring{\phi}=-\varkappa^2\mathring{\phi}.
\end{equation}
This differential equation has a general solution of the form
$\mathring{\phi}=Ae^{\varkappa t}+Be^{-\varkappa t}$. 
Since we want that our solution has similar behavior to the so-called
flat de Sitter solution in general relativity (see for example
\cite[Chap.~V]{Choquet-Bruhat_09
}), namely, that $\mathring{\phi}$ and $\frac{d}{dt}\mathring{\phi}$ are
positive, we chose
\begin{equation}
  \label{eq:section1-intro:2B}
  \mathring{\phi}(t)= \mathord{e^{\varkappa t}}
\end{equation}
as the background solution. 
Considering also the part $Be^{-\varkappa t}$ would complicate the
analysis but should not change the global behavior of the solutions,
that is why we are neglecting this term.

We now study small deviations from the background solution
$\mathring \phi$. 
So we make the following Ansatz
\begin{equation}
 \label{eq:phi}
 \phi=\mathring{\phi}+\Psi=\mathord{e^{\varkappa t}} +\Psi,
\end{equation}
where $\Psi$ denotes the deviation from the background. 
Then $\Psi$ satisfies the following equation
\begin{equation}
  \square \phi=\square(\mathord{e^{\varkappa t}}+\Psi)=-\varkappa^2 \mathord{e^{\varkappa t}}+\square 
  \Psi=-\frac{1}{6} T(\mathord{e^{\varkappa t}} 
+\Psi)^3-\varkappa^2(\mathord{e^{\varkappa t}} +\Psi).
\end{equation}
Thus $\Psi$ satisfies the initial value  problem
\begin{equation}
\label{eq:psi}
 \left\{
 \begin{array}{l}
 \square \Psi =-\frac{1}{6}T(\mathord{e^{\varkappa 
t}}+\Psi)^3-\varkappa^2\Psi\\
 \Psi(0,x)=\Psi_0(x), \ 
   \partial_t\Psi(0,x)=\Psi_1(x)
 \end{array}\right..
\end{equation} 
Our goal is:  
\begin{enumerate}[label=\alph*.] 
  \item \label{item:section1-intro:2}
  To show global existence of classical solutions for equation \eqref{eq:psi}
  demanding a small source term $T$ and small initial data.
  \item
  \label{item:section1-intro:1} To show that for large $t$, the metric
  $\phi^2\eta_{\alpha\beta}$ approaches  asymptotically  the background metric
  $e^{2\varkappa t}\eta_{\alpha\beta}$, in the following sense, 
  \begin{equation}
    \label{eq:section1-intro:3B}
    \lim_{t\to\infty}\frac{\phi(t,x)}{\mathring      \phi(t)}=\lim_{t\to\infty}\frac{e^{\varkappa t}+\Psi(t,x)}{e^{\varkappa t}}\approx 1.
  \end{equation}

\end{enumerate}

Note that if $\Psi$ is small, then
$(\mathord{e^{\varkappa t}} +\Psi)^3\sim e^{3\varkappa t}$, and this term
growths very rapidly and might prevent that the solution exists for
all time. 
So in order to achieve the desired asymptotic behavior of $\Psi$,
expressed by equation \eqref{eq:psi}, we multiply $\phi$ by
$\mathord{e^{-\varkappa t}}$, then from equality \eqref{eq:phi} we
conclude that
$\mathord{e^{-\varkappa t}} \phi=1+\mathord{e^{-\varkappa t}}\Psi $, and
therefore we set
\begin{equation}
  \Omega\overset{\mbox{\tiny{def}}}{=} \mathord{e^{-\varkappa t}} \Psi.
\end{equation}
The resulting equation for $\Omega$ takes the form
\begin{align*}
  \partial_t \Omega &= \partial_t (e^{-\varkappa t} \Psi)=e^{-\varkappa t}
              \partial_t\Psi-\varkappa e^{-\varkappa t} \Psi=e^{-\varkappa t} \partial_t\Psi-\varkappa \Omega,\\
  \partial_t^2 \Omega &=  e^{-\varkappa t} 
               \partial_t^2\Psi-2\varkappa e^{-\varkappa t}\partial_t\Psi+\varkappa^2e^{-\varkappa t}\Psi= 
               e^{-\varkappa t}\partial_t^2\Psi-2\varkappa(\partial_t \Omega+\varkappa 
               \Omega)+\varkappa^2\Omega\\
            & = e^{-\varkappa t}\partial_t^2\Psi-2\varkappa \partial_t 
              \Omega -\varkappa^2\Omega,
\end{align*}
or
\begin{equation}
  -e^{-\varkappa t}\partial_t^2\Psi=-\partial_t^2\Omega-2\varkappa \partial_t 
\Omega -\varkappa^2\Omega.
\end{equation} 
Thus we have obtained
\begin{equation}
 e^{-\varkappa t} \square \Psi=\square \Omega-2\varkappa\partial_t 
\Omega-\varkappa^2\Omega=-\frac{1}{6} \widehat Te^{-\varkappa t}(e^{\varkappa 
t}+\Psi)^3-\varkappa^2\Omega,
\end{equation} 
or
\begin{equation}
  \label{eq:Omega}
      \square \Omega-2\varkappa\partial_t \Omega=-\frac{1}{6}  
T(t,x)e^{2\varkappa t}(1+\Omega)^3.
\end{equation}

We wish to show the existence of global classical solutions for system
\eqref{eq:Omega} demanding a small source term $T(t,x)$.  
On the one hand, the term $e^{2\varkappa t }$ seems to hamper the proof
of the desired global existence, but on the other hand, we have
obtained a good dissipative term of the form $-2\varkappa\partial_t\Omega$.
This is why we perform the transformation
$\widetilde T= \widetilde g_{\alpha\beta}\widetilde T^{\alpha\beta}=
e^{3\kappa t} T$, which also implies that the right-hand of the wave
equation \eqref{eq:Omega} takes the form
$-\frac{1}{6}e^{-\varkappa t} \widetilde T(1+\Omega)^3$. 
If $\Omega$ remains bounded, then the right-hand side will tend to zero. 
That is why we finally consider the following system
\begin{subequations}
  \begin{align}
    \label{eq:section2-field:1}
    - &\partial_t^2 \Omega-2\varkappa \partial_t \Omega +\Delta 
        \Omega=-\mathord{e^{-\varkappa t}} a(t,x)(1+\Omega)^3\\ 
    \label{eq:section2-field:2}
      & (\Omega(0,x),\partial_t \Omega(0,x))=(f(x),g(x)),
  \end{align}
\end{subequations}
where we have denoted $\frac{1}{6}\widetilde T$ by $a(t,x)$.

\begin{rem}[The scaling and the Euler equations]
  The above scaling of the trace of the energy-momentum tensor will
  change the Euler equations. 
  That is why this scaling has to be taken into account for the
  coupled Euler-Nordström system, which we want to treat in a
  forthcoming paper. 
  Moreover, it turns out that we also need to scale the metric and the
  velocity as follows:
  $\widetilde g_{\alpha\beta}=e^{-2\kappa t} g_{\alpha\beta}=e^{-2\kappa
    t}\phi^2 \eta_{\alpha\beta}$, and the
  $\widetilde{u}^\alpha =\mathord{e^{\varkappa t}} u^\alpha$, which is
  compatible with the scaling $\widetilde T= e^{3\kappa t} T$.
 
\end{rem}

In what follows we will not consider the Euler-Nordström system but
instead consider the fluid as a given source of the field equations,
and therefore we will consider the right-hand side of equation
\eqref{eq:Euler-Nordstrom:5} as a given function of $(t,x)$.
A similar setting was considered by H.~Friedrich for the Einstein
vacuum equations with positive cosmological constant, in which he
proved global existence of classical solutions for small initial data
\cite{Friedrich:1986
}.

We point out that we require the deviation
$e^{\varkappa t}\Omega=\Psi$ to be spatially periodic, and that is why
we study the Cauchy problem
\eqref{eq:section2-field:1}--\eqref{eq:section2-field:2} in the Sobolev
spaces $H^m(\setT^3)$. 
We shall decompose this space into two orthogonal components, namely,
$H^m(\setT^3)=\setR \oplus \mathbullet H^m(\setT^3)$, where the second
component consists of all functions with zero mean over the torus
$\setT^3$. 
The reason for this decomposition is that the homogeneous space
$\mathbullet H^{m}(\setT^3)$ possesses some convenient features for our energy
estimates and seems best suited for our setting.
However, there is a technical difficulty in using these spaces, namely
the presence of the nonlinear term $(1+\Omega)^3$, that cannot belong to the
homogeneous spaces. 
We solve this problem by performing a projection of our variables into
a part that belongs to these spaces, and another part that satisfies an
ordinary differential equation, we refer to section
\ref{sec:math-prel} for details.
As we will see, in section \ref{sec:math-prel}, these spaces posses
some nice features, such as Proposition \ref{prop:9}, that simplify
the energy estimates which we shall use for proving our results.

Having set up the problem, we outline the structure of our paper and
summarize our main results.

In section \ref{sec:math-prel}, we introduce the necessary
mathematical tools, such as homogeneous and non-homogeneous Sobolev
spaces on the torus $\setT^3$.
Using Fourier series, in section \ref{sec:semi-linear-wave}, we
obtain, for small initial data and a small source term, global
existence and uniqueness  of these solutions in the
$ H^m(\setT^3)$ spaces (see Theorem \ref{thm:1}).

We then turn, in Section \ref{symm-hyper-const}, to the theory of
symmetric hyperbolic systems. 
We write the wave equation in a slightly unorthodox way as a symmetric
hyperbolic system (see system \eqref{eq:symm:1}) and then prove global
existence, uniqueness and asymptotic decay for a small source term,
but not necessarily small initial data,
(see~Theorem~\ref{thr:section5-sh:1}). 
The reason we consider the semi-linear wave equation
\eqref{eq:section2-field:1} in the framework of the theory of
symmetric hyperbolic systems is that in the future we want to consider
the coupled Euler-Nordström system, and we know already the Euler
equations can be cast into that form (see
\cite{BK8
}).
Finally, in the Section \ref{sec:blow-solut-large}, we show that if
the source term is not small, then the corresponding solutions blow up
in finite time. 
It turns out, however, that for the proof of the blow up result we
need that $(1+\Omega(t,x))\geq 0$, which seems natural if the initial data are
positive.  
However, for that being true, it is not sufficient to only assume the
initial data to be positive; additional conditions are needed that also
result in a more elaborated proof.
That has been taken care of in the last section.

\biblio

%% file: section2-preliminaries.tex
\section{Mathematical Preliminaries}
\label{sec:math-prel}

\subsection{Sobolev spaces on the torus $\setT^3$}
\label{sec:non-homog-sobol}

We consider the solutions on the torus $\mathbb T^3$ using 
Sobolev spaces $H^m$ where $m$ is a nonnegative integer (see
e.g.~\cite[Chap~3.1]{taylor97
},
\cite[Chap.~5.10]{robinson2001
}).
It is natural to represent  functions  on the torus by Fourier series
and their norms by Fourier coefficients.  
For a function $f$, its Fourier series is given by
\begin{equation}
 \label{eq:fs:1}
 f(x)=\sum_{ k\in \setZ^3}  \widehat f_{ k}e^{ix\cdot  k},
\end{equation} 
where 
\begin{equation}
 \label{eq:fs:2}
 \widehat f_{ k}=\frac{1}{(2\pi)^3}\int_{\mathbb T^3}f(x)e^{-ix\cdot  k } dx,
\end{equation}
$x\cdot k=x_1k_1+x_2k_2+x_3k_3$, $x\in \setT^3$, and $ k\in \setZ^3$.

The $H^m$ norm is given by
\begin{equation}
 \label{eq:norm:1}
 \|f\|^2_{H^m(\mathbb T^3)}=\|f\|^2_{H^m}= |\widehat f_0|^2+ 
\sum_{ k\in\setZ^3}| k|^{2m}|\widehat f_{ k}|^2.
\end{equation}

The homogeneous Sobolev spaces $\mathbullet H^m$ are defined by the semi--norm
\begin{equation}
 \label{eq:norm:2}
 \|f\|^2_{\mathbullet H^m(\mathbb T^3)}=\|f\|^2_{\mathbullet H^m}=  
\sum_{ k\in\setZ^3}| k|^{2m}|\widehat f_{ k}|^2.
\end{equation}
We decompose the Sobolev space $H^m$ into two orthogonal components
\begin{equation}
 \label{eq:norm:3}
 H^m=\setR \oplus \mathbullet H^m.
\end{equation}
A function $f\in H^m$ belongs to $\mathbullet H^m$ if and only if it has a zero 
mean, that is, 
\begin{equation}
 \label{eq:fs:3}
 \frac{1}{(2\pi)^3}\int_{\mathbb T^3}f(x)dx=0.
\end{equation} 

By Parseval's identity, 
\begin{equation*}
 \|f\|^2_{L^2(\mathbb T^3)}=\sum_{k\in\setZ^3}|\widehat f_k|^2
\end{equation*}
and since $\widehat{(\partial^\alpha f)}_k =k^\alpha \widehat{f}_k$, the 
following equivalent holds
\begin{equation*}
 \|f\|^2_{ H^m( \setT^3)}\simeq \|f\|^2_{L^2(\mathbb T^3)}+ 
\sum_{|\alpha |= m} \|\partial^\alpha 
f\|_{L^2(\setT^3)}^2.
\end{equation*}

We also introduce an inner-product,  for two vector 
valid real functions $U$ and $V$, we set
\begin{equation}
\label{eq:norm:4}
 \langle U, V\rangle_m=\widehat U_0\cdot \widehat V_0+ \sum_{k\in 
\setZ^3}|k|^{2m}\left(\widehat U_k\cdot \overline{\widehat V_k}\right).
\end{equation}

The following proposition which is a certain version of Wirtinger's
inequality \cite{Dym_P.Mckean_85
} is a simple consequence of the representation of the homogeneous
norm \eqref{eq:norm:2}. 
\begin{prop}[Estimate for the gradient]
  \label{prop:9}
  Let $ \partial_x u=(\partial_1 u,\partial_2u,\partial_3 u)^{\intercal}$ and
  $u\in \mathbullet H^{m+1}$ and $m\geq0$ be an integer. 
  Then the following holds
  \begin{equation*}
    \label{eq:norm:36}
    \left\|{u}\right\|_{\mathbullet H^{m+1}}=\|\partial_x u\|_{\mathbullet       H^{m}}.
  \end{equation*}

\end{prop}
 
\begin{proof}
  Since $|(\widehat{\partial_x u})_k|^2=|k|^2\ |\widehat u_k|^2$, we obtain by the
  representation \eqref{eq:norm:3} of the norm that
  \begin{equation}
    \|\partial_x u\|_{\mathbullet H^{m}}^2=\sum_{0\neq k\in\setZ^3} 
    |k|^{2m} 
    |(\widehat{\partial_x u})_k|^2 =\sum_{0\neq k\in\setZ^3} 
    |k|^{2m} 
    |k|^2 |\widehat u_k|^2=\left\|u\right\|_{\mathbullet H^{m+1}}^2
  \end{equation}
which proves the proposition.  
\end{proof}

 \subsection{Calculus in the Sobolev spaces on the torus $\setT^3$}
\label{sec:calculus-sobol}
 
We recall that the known properties in Sobolev spaces, defined over
$\setR^n$ such as  multiplication, embedding and Moser type
estimates, hold also for Sobolev spaces defined over the torus
$\setT^n$, see e.g.~\cite[Chap.~13]{taylor97c
}.

\begin{prop}[A Nonlinear estimate]
\label{prop:1}
 Let $m>\frac{3}{2}$ and  $a  \in  H^{m}$, then there is a universal constant 
$C(A)$, depending just on the constants of multiplications and embedding, such 
that 
 \begin{equation}
   \label{eq:estimate:3}
   \left\|a(1+u)^3 \right\|_{H^{m}}, \left\|a(1+u)^3 \right\|_{L^\infty}\leq  
C(A) \|a\|_{H^{m}} 
 \end{equation}
for all $u\in H^m$ with $\|u\|_{H^m}\leq A$.
 
\end{prop}

\begin{proof}
  By the multiplication property (\cite[Proposition 3.7,~Chap.~13]{taylor97c
  }), there is a constant $C$ such that 
  \begin{equation}
    \label{eq:section3-fourier:1}
    \begin{split}
    \left\|a(1+u)^3\right\|_{H^{m}} & \leq C
    \left\|a\right\|_{H^{m}} 
    \left\|(1+u)^3\right\|_{H^{m}}  
    \leq C \left\|a\right\|_{H^{m}} 
    \left\|(1+u)\right\|_{H^{m}}^3  \\ & \leq C \left\|a\right\|_{H^{m}} 
    \left(1+\|u\|_{H^m}\right)^3 \leq C \left\|a\right\|_{H^{m}} 
    \left(1+A\right)^3.
 \end{split}
 \end{equation}
 Using the embedding  $\|u\|_{L^\infty}\leq C \|u\|_{H^m}$, we  see that 
\eqref{eq:estimate:3} holds.
 
\end{proof}

\subsection{Estimate of symmetric hyperbolic system}

We shall also need the following property of solution to semi-linear
symmetric hyperbolic systems.
Consider a symmetric hyperbolic system
\begin{equation}
\label{eq:sh:1}
 \partial_t U=\sum_{j=1}^3A^j(t,x)\partial_j U +F(t,x,U),
\end{equation} 
where the matrices $A^j(t,x)$ are symmetric and $F(t,x,U)$ is a smooth function 
of $U $. The next proposition provides a uniform modulus of continuity for the 
difference $U(t,\cdot)-U_0(\cdot) $ in the $H^m $ norm.

\begin{prop}[Modulus of continuity]
  \label{prop:6}
  Let $m>\frac{5}{2}$, $A^j \in L^\infty([0,T];H^m )$ for some positive
  $T$ and $F(t,x,0)\in L^\infty([0,T];H^m )$. 
  Assume that $U(t)\in C([0,T];H^m )\cap C^1([0,T]);H^{m-1})$ is the
  solution to system \eqref{eq:sh:1} with initial data $U_0\in H^m$,
  then there is a constant $C(\|U_0\|_{H^m})$ such that
  \begin{equation}
    \|U(t)-U_0\|_{ H^{m-1}}\leq 
    C({\left\|{U_0}\right\|_{H^m}})t^{\frac{1}{m}}\quad \text{for}\ \  0<t<T. 
  \end{equation}
\end{prop}

\begin{rem}

We know that from the existence theory for quasilinear symmetric hyperbolic 
system that the solution $U$ belongs to a certain ball around $U_0$ in the 
$H^m$ 
space (see e.g. \cite{KATO}, \cite{rauch12:_hyper}). So we may assume that 
$\|U(t)\|_{H^m}\leq 
\|U_0\|_{H^m}+R$ for some positive $R$ and $t\in[0,T]$. The same phenomena 
appears also for quasi-linear wave equations (see e.g. 
\cite[Theorem 6.4.11]{Hormander_1997}).
 
\end{rem}

\begin{proof}
  Let $ t<T$, then
  \begin{equation}
    U(t,x)-U_0(x)=\int_{0}^t\partial_tU(\tau,x) d\tau.
  \end{equation}
  By the Cauchy Schwarz inequality, it follows that 
  \begin{equation}
    \Big| U(t,x)-U_0(x)\Big|^2\leq t\int_{0}^t|\partial_tU(\tau,x)|^2 
    d\tau.
  \end{equation}
  Hence we conclude that 
  \begin{equation}
    \label{eq:mod:4}
    \| U(t,x)-U_0(x)\|_{L^2}^2 \leq  t\int_{0}^t\int 
    |\partial_tU(\tau,x)|^2dx d  \tau \leq
    t^2\left\|\partial_t U\right\|_{L^\infty([0,T];L^2)}^2.
  \end{equation}
  Since $U$ satisfies system \eqref{eq:sh:1}, we obtain
  \begin{equation}
    \label{eq:sh:2}
    \begin{split}
  \|\partial_t U(t)\|_{L^2} &\leq     \|\partial_t U(t)\|_{H^{m-1}}  \leq 
\sum_{j=1}^3 \|A^j(t,x)\partial_j U(t)\|_{H^{m-1}}+
      \|F(t,x,U(t))\|_{H^{m-1}}\\ &  \leq  C\sum_{j=1}^3 \|A_j(t,\cdot)\|_{H^m} 
\|U(t)\|_{H^{m}} +
     C(\|U(t)\|_{L^\infty}) \|U(t)\|_{H^m}+ \|F(t,\cdot,0)\|_{H^m}.
    \end{split}
  \end{equation}
  Here we used the multiplication property and Moser third estimate,
  see e.g.
  \cite[Theorem~6.4.1]{rauch12:_hyper
  }, \cite[Proposition 3.9,~Chap.~13]{taylor97}.
  Thus it follows from the remark that
 \begin{equation}
  \sup_{[0,T]}\|\partial_t U(t)\|_{L^2}\leq C(\|U_0\|_{H^m}).
 \end{equation} 

 We now apply the intermediate estimate
 $\|u\|_{H^r}\leq \|u\|_{H^{m}}^{m-\frac{r}{m}}\|u\|_{L^2}^{\frac{r}{m}}$ for
 $0<r< m$, (see e.g.~\cite[Prop.~1.52]{Bahouri_2011}) and inequality
 \eqref{eq:mod:4}, then
   \begin{equation}
    \label{eq:estimate:9}
    \begin{split}
      \|U(t)-U_0\|_{ H^{m-1}}\leq \|U(t)-U_0\|_{ H^{m}}^{\frac{(m-1)}{m}}
      \|U(t)-U_0\|_{L^2}^{\frac{1}{m}} \leq \|U(t)-U_0\|_{H^m}^{frac{m-1}{m}}
      C_0(\|U_0\|_{H^m})t^{\frac{1}{m}}.
    \end{split}
  \end{equation} 
  Since $\|U(t)\|_{H^m}\leq \|U_0\|_{H^m}+R$, that completes the proof.
\end{proof}

\subsection{Gronwall inequality}
We shall use the following version of Gronwall's inequality (see 
e.~g.~\cite{Bahouri_2011}).
\begin{lem}[Gronwall's inequality]
  \label{lem:Gronwall}
  Let $g$ be a $C^1$ function, $f$, $F$, and $A$ continuous function in
  the interval $[t_0, T]$. 
  Suppose that for $t\in [t_0,T]$ $g$ obeys
  \begin{equation}
\label{eq:section2-preliminaries:1}
    \frac{1}{2}\frac{d}{dt}g^2(t)\leq A(t)g^2(t)+f(t)g(t).
  \end{equation}
  Then for $t\in [t_0,T]$ we have
  \begin{equation}
\label{eq:section2-preliminaries:2}
    g(t)\leq e^{\int_{t_0}^t A(\tau)d\tau} g(t_0)+\int_{t_0}^t 
    e^{\int_{\tau}^tA(s)ds} f(\tau)d\tau.
  \end{equation}

\end{lem}

\biblio

%% file: section3-fourier.tex
\section{The Cauchy problem for a semi-linear wave equation using
  Fourier series}
\label{sec:semi-linear-wave}

In the following section we shall investigate the Cauchy problem
\eqref{eq:section2-field:1}--\eqref{eq:section2-field:2}, however, for
convenience, we multiply the wave equation by $-1$ and denote the
unknown by $u$ instead of $\Omega$, which results in the following
semi-linear wave equation
\begin{subequations}
  \begin{align}
    \label{eq:wave:1}
    &\partial_t^2 u+2\varkappa \partial_t u -\Delta u=\mathord{e^{-\varkappa t}} a(t,x)(1+u)^3\\ 
    \label{eq:wave:2}
    & (u(0,x),\partial_t u(0,x))=(f(x),g(x)).
  \end{align}
\end{subequations}

Here $\varkappa$ is a positive constant, while $a(t,x)$ is a smooth function
as we discussed in Section \ref{sec:fields-equations}.

We are interested in proving the global existence of classical solutions
to the Cauchy problem \eqref{eq:wave:1}--\eqref{eq:wave:2} for small
initial data and $a(t,x)$.
We also note, that the Cauchy problem
\eqref{eq:wave:1}--\eqref{eq:wave:2} has some similarities with the
Cauchy problem of the damped semi-linear wave equation
\begin{subequations}
  \begin{align}
    \label{eq:section3-fourier:7}
    & \partial_t^2 u+2\varkappa \partial_t u -\Delta u=|u|^p, \qquad (t,x)\in \setR_+\times \setR^3\\
    \label{eq:section3-fourier:6}
    & (u(0,x),\partial_t u(0,x))=(f(x),g(x)).
  \end{align}
\end{subequations}
for which it is known that for $2\leq p\leq 3$ there exist global solutions
for small initial data, for further details we refer to
\cite{Ebert_Reissig_18
}, \cite{Todorova_Yordanov_01
} and the references therein.
However, we did not find in the literature any results concerning the
Cauchy problem
\eqref{eq:section3-fourier:7}--\eqref{eq:section3-fourier:6} on the
torus.
There is however another difference, between these two Cauchy problems
\eqref{eq:section3-fourier:7}--\eqref{eq:section3-fourier:6} and
\eqref{eq:wave:1}--\eqref{eq:wave:2}.
In equation \eqref{eq:wave:1}, the function $a(t,x)$ is essential, in the sense
that global existence depends on the smallness of this function, while
the structure of the nonlinear term is of less importance.

\subsection{Local existence}
\label{sec:local-existence}
Before we are going to present our results concerning the global
existence of classical solutions we shall discuss the question of
local existence and uniqueness of the initial value problem
\eqref{eq:wave:1}--\eqref{eq:wave:2}. 
There are well known local existence and uniqueness theorems for
quasilinear wave equations of the form
\begin{equation}
\label{eq:section3-fourier2}
 g^{\alpha\beta}(u,u')\partial_\alpha\partial_\beta u=F(u),
\end{equation}
where $g^{\alpha\beta}(u,u')$ has a Lorentzian signature and
$u'=\partial_\alpha u$, $\alpha=0,1,2,3$, see for example \cite[Theorem
6.4.11]{Hormander_1997}, \cite[Theorem 4.1]{Sogge_95} and
\cite[Theorem 5.1]{Shatah-Struwe-98}.
These references treat the initial value problems in the Sobolev space
$H^m(\setR^3)$ and under the condition $F(0)=0$. 
We consider solutions of equation \eqref{eq:wave:1} that belong to
Sobolev spaces on the torus $\setT^3$, $H^m(\setT^3)$, and we observe that the
right hand side of \eqref{eq:section3-fourier2} does not satisfy the
condition $F(0)=0$. 
Nevertheless, the above existence results can be applied to the Cauchy
problem \eqref{eq:wave:1}--\eqref{eq:wave:2} because of the following
reasons.

\begin{enumerate}[label=\arabic*.]
  \item The energy estimates are an indispensable tool for proving
  local existence for the linearized equation. 
  The energy estimates rely on the formula for integration by parts
  $\int u\partial_{x_j} v dx=-\int \partial_{x_j} u v dx$, which holds for periodic functions
  and rapidly decreasing functions in $\setR^n$. 
  That is why the energy estimates in the above references hold in the
  Sobolev spaces $H^m(\setT^3)$ as well.  

  \item Moser type inequality, the second important tool, states that
  $\|F(u)-F(0)\|_{H^m}\leq C \|u\|_{H^m}$ for a sufficiently smooth function
  $F$. 
  This nonlinear estimate is valid both for $u\in H^m(\setR^3)$ and for
  $u\in H^m(\setT^3)$. 
  In the case the equations are considered on the $\setR^3$, the
  requirement $F(0)=0$ is needed since the constant function does not
  belong to Sobolev space $H^m(\setR^3)$. 
  However, the situation is different on the torus. 
  Here obviously the constant function belongs to the space
  $H^m(\setT^3)$. 

\end{enumerate}

So we conclude that with some minor modifications of \cite[Theorem
6.4.11]{Hormander_1997}, the following result on local existence and
uniqueness.

\begin{thm}[Local existence]
  \label{thr:section3-fourier:1}
  Let $m>\frac{5}{2}$,
  $a(t,\cdot)\in L^\infty([0,\infty); H^m(\setT^3)$,
  $ f\in H^{m+1}(\setT^3)$ and $ g\in H^{m}(\setT^3)$, then there exists a positive
  $T$ and a unique solution $u$ to the Cauchy problem
  \eqref{eq:wave:1}--\eqref{eq:wave:2} such that
  \begin{equation*}
    u\in L^\infty([0,T]; H^{m+1}(\setT^3)\cap    C^{0,1}([0,T]; H^{m}(\setT^3), 
  \end{equation*}
  where $C^{0,1}$ is a Lipschitz continuous function.
\end{thm}

\subsection{Global existence}
\label{subsec:homogeneous}
Once that the existence and uniqueness of local solutions have been
established (by theorem \ref{thr:section3-fourier:1}), we turn now to
the question of whether global solutions to the Cauchy problem
\eqref{eq:wave:1}--\eqref{eq:wave:2} exist.
Again, the energy estimates are the main tool for treating this
problem.
Those energy estimates are different from the one that has been used
for local existence. 
Our method consists in expanding the solution into Fourier series,
which allows us to solve the corresponding ordinary differential
equations for the Fourier's coefficients, and use the norm
\eqref{eq:norm:1} to derive the desired estimates.   

Our main result in  this section is the following theorem.

\begin{thm}[Global existence of classical solutions for small data]
  \label{thm:1}
  Let $m> \frac{5}{2}$, $f\in  H^{m+1}$,
  $g\in  H^m$, and $a\in C([0,\infty); H^{m})$. 
  Then there is a suitable constant $\varepsilon$ such that if the following holds
  \begin{equation*}
    \|f\|_{ H^{m+1}}, \|g\|_{ H^{m}}, 
    \sup_{[0,\infty)}\|a(t,\cdot)\|_{ H^{m}}<\epsilon, 
  \end{equation*}
  then the Cauchy problem \eqref{eq:wave:1}--\eqref{eq:wave:2} has a
  unique global solution of the form
  \begin{equation}
    \label{eq:section3-fourier:4}
    u\in C([0,\infty); H^{m+1}).
  \end{equation}
  Moreover, there exists  a positive constant $C_1$ such that
  \begin{equation}
    \label{eq:limit-zero}
    \|u(t,\cdot)\|_{ H^{m+1}}\leq {C_1}.
  \end{equation}
 \end{thm}

\subsection{Proof of  Theorem \ref{thm:1}}
\label{sec:outline-proof}

The main points and ideas of the proof of Theorem \ref{thm:1} can be
described as follows:
\begin{enumerate}[label=\arabic*.]
  \item Obtain an energy estimate for the linearized equation.
  \item Use the Banach fixed pointed theorem for the linearized
  equation. 
\end{enumerate}

We start  with the energy estimates for the linearized system
of equation
\eqref{eq:wave:1}. 
 For any function
$v\in  H^{m}$ we set
\begin{equation}
  \label{eq:section3-fourier:9}
  F(t,x)\eqdef a(t,x)\left( 1+v \right)^{3}
\end{equation}
   and consider the following linear initial value problem
\begin{subequations}
  \begin{align}
    \label{eq:wave:7}
    &\partial_t^2 u+2\varkappa \partial_t u -\Delta u= \mathord{e^{-\varkappa t}} 
      F(t,x)\\ 
    \label{eq:wave:8}
    & (u(0,x),\partial_t u(0,x))=(f(x)),g(x)).
  \end{align}
\end{subequations}

We now consider  the Fourier coefficients 
\begin{equation*}
 \widehat u_k(t)=\frac{1}{(2\pi)^3}\int_{\setT^3} u(t,x)dx, \quad k\in \setZ^3,
\end{equation*}
and the coefficients of the other data of the Cauchy problem
\eqref{eq:wave:7}--\eqref{eq:wave:8} as well. 
We obtain, by this procedure, for each $ k\in \setZ^3$ an ordinary
differential equation
\begin{subequations}
  \begin{align}
    \label{eq:fs:4}
    & \widehat u^{\prime\prime}_{ k}(t)+2\varkappa \widehat u'_{ k}(t)+| 
      k|^2\widehat u_{ 
      k}(t)=\mathord{e^{-\varkappa 
      t}} \widehat F_{ k}(t)\\
    \label{eq:fs:5}
    & \widehat u_{ k}(0)=\widehat f_{ k}, \quad \widehat u'_{ k}(0)=\widehat 
      g_{ k}.
  \end{align}
\end{subequations}
We can solve \eqref{eq:fs:4}--\eqref{eq:fs:5} explicitly, however,
since the structure of the solutions depends on $\varkappa$, and in
order to work with similar formulas for all $k\neq 0$, we restrict
$\varkappa$ to the interval $(0,1)$.
We present the energy estimates in the following proposition
\begin{prop}[Energy estimate for the linearized wave equation]
  \label{prop:3.2}
  Let $ m\geq 0$ and $0<\varkappa<1$, and assume
  $F\in C([0,\infty); H^m)$, $f\in  H^{m+1}$, and
  $g \in H^m$. 
  Then there exists a unique solution
  $u\in C([0,\infty); H^{m+1})$ to equation \eqref{eq:wave:7} with
  initial data \eqref{eq:wave:8}, and moreover, it obeys
  \begin{equation}
    \label{eq:estimate-linear}
    \begin{split}
      \|u(t,\cdot)\|^2_{H^{m+1}} & \leq {2e^{-2\varkappa
          t}}\left\{(1+2\varkappa^2)(1+t^2)\|f\|^2_{\mathbullet H^{m+1}}+(2(1+t^2)
        \|g\|^2_{\mathbullet H^{m}}+t(1+t^2)\int_0^t
        \|F(\tau,\cdot)\|^2_{\mathbullet H^{m}} d\tau\right\} \\ & + \widehat
      f_{0}^2+\widehat g^2_{0}\left(\frac{1-e^{-2\varkappa
            t}}{2\varkappa}\right)^2+\frac{1}{4} \sup_{\tau\in [0,t]}|\widehat
      F_{ 0}(\tau)|^2\left( \frac{1-e^{-\varkappa t}}{\varkappa}\right)^4.
    \end{split}
  \end{equation}

\end{prop}

\begin{proof}
  For each $k\neq 0$ the solution of the initial value problem of the
  ordinary differential equations \eqref{eq:fs:4}--\eqref{eq:fs:5} is
  then given by
  \begin{equation}
    \label{eq:fs:7}
    \begin{split}
      \widehat u_{ k}(t) & =\mathord{e^{-\varkappa t}} \left\{ \widehat
        f_{ k} \cos\left(\sqrt{|k|^2-\varkappa^2}\, t\right)+\frac{\widehat
          g_{ k}+\varkappa \widehat f_{ k}}{\sqrt{|k|^2-\varkappa^2}}\,
        \sin\left(\sqrt{|k|^2-\varkappa^2}\, t\right)\right\} \\ &
      +\frac{1}{\sqrt{|k|^2-\varkappa^2}}\int_0^t
      e^{-\varkappa(t-\tau)}\sin\left(\sqrt{|k|^2-\varkappa^2}\,(t-\tau)\right)e^{
        -\varkappa \tau}\widehat       F_{ k}(\tau)d\tau\\
      & = \mathord{e^{-\varkappa t}} \left\{ \widehat f_{ k}
        \cos\left(\sqrt{|k|^2-\varkappa^2}\,t\right)+\frac{\widehat g_{
            k}+\varkappa \widehat f_{
            k}}{\sqrt{|k|^2-\varkappa^2}}\sin\left(\sqrt{|k|^2-\varkappa^2}
          t\right)\right. 
      \\ & +\left.        
        \frac{1}{\sqrt{|k|^2-\varkappa^2}}\int_0^t\sin\left(\sqrt{|k|^2-\varkappa^2}
          (t-\tau)\right)\widehat F_{ k}(\tau)d\tau \right\},
    \end{split}
  \end{equation}
  and for $ k=0$,
  \begin{equation}
    \label{eq:u-0}
    \widehat u_{ 0}(t)=\widehat f_{ 0}+\widehat 
    g_{0}\left(\frac{1-e^{-2\varkappa        
          t}}{2\varkappa}\right)+\frac{1}{2\varkappa}\int_0^t\left(1-e^{
        -2\varkappa(t-\tau) }
    \right)e^{-\varkappa \tau} \widehat F_{ 0}(\tau) d\tau.
  \end{equation}

  We shall now estimate $\|u\|_{H^{m+1}}^2$ by the formula
  \eqref{eq:norm:1}. 
  For $k\neq 0$, we conclude from equality \eqref{eq:fs:7} and the trivial 
inequality $(a+b+c)^2\leq 2(a^2+b^2+c^2)$ that
  \begin{equation*}
    \begin{split}
      & |k|^{2(m+1)}|\widehat u_{ k}(t)|^2 \leq 2e^{-2\varkappa
        t}|k|^{2(m+1)}\left\{|\widehat
        f_k|^2\left(\cos\left(\sqrt{|k|^2-\varkappa^2}\, t\right)\right)^2
      \right. 
      \\ + & \left. 
        \left|\widehat g_k+\varkappa \widehat
          f_k\right|^2\left(\frac{\sin\left(\sqrt{|k|^2-\varkappa^2}\,
              t\right)}{\sqrt{|k|^2-\varkappa^2}}\right)^2+ \left(\int_0^t
          \frac{\sin\left(\sqrt{|k|^2-\varkappa^2}\, (t
              -\tau)\right)}{\sqrt{|k|^2-\varkappa^2}}
          \widehat{F_k}(\tau)d\tau\right)^2 \right\}\\ & = I_k +II_k+III_k.
    \end{split}
  \end{equation*}

  The first one is easy to estimate, and we obtain that
  \begin{equation}
    \label{eq:e-1}
    I_k\leq 2e^{-2\varkappa t}|k|^{2(m+1)}|\widehat f_k|^2.
  \end{equation}
  For the second and third term we use the inequality
  $\sqrt{1+\xi^2}|\sin(\xi t)|\leq \xi\sqrt{1+t^2}$, with
  $\xi=\sqrt{|k|^2-\varkappa^2}$, that implies
  \begin{equation}
    \frac{\sin\left(\sqrt{|k|^2-\varkappa^2}\,      
t\right)}{\sqrt{|k|^2-\varkappa^2}}\leq\sqrt{\frac{1+t^2}{1+\xi^2}}=\sqrt{\frac{
        1+t^2}{|k|^2+1-\varkappa^2}}\leq \frac{\sqrt{1+t^2}}{|k|}.
  \end{equation}
  Hence
  \begin{equation}
    \label{eq:e-2}
    II_k\leq 2 e^{-2\varkappa t}|k|^{2m}|\widehat g_k +\varkappa \widehat 
    f_k|^2(1+t^2)\leq 2 e^{-2\varkappa t}|k|^{2m}\left(2|\widehat g_k|^2 
      +2\varkappa^2 |\widehat 
      f_k|^2\right)(1+t^2)
  \end{equation}
  and
  \begin{equation}
    \label{eq:e-3}
    III_k\leq 2 e^{-2\varkappa t}|k|^{2m} t\int_0^t {(1+(t-\tau)^2}|\widehat 
    F_k(\tau)|^2d\tau\leq 2 e^{-2\varkappa t}|k|^{2m}|
    t(1+t^2)\int_0^t |\widehat 
    F_k(\tau)|^2d\tau.
  \end{equation}

  We now turn to the zero's term \eqref{eq:u-0}, and we start with the
  integral term of \eqref{eq:u-0},
  \begin{equation}
    \label{eq:u-0:2}
    \bigg|\frac{1}{2\varkappa}\int_0^t\left(1-e^{-2\varkappa(t-\tau)}\right) 
    e^{-\varkappa \tau} \widehat F_{ 0}(\tau) d\tau \bigg|\leq 
    \frac{1}{2}\left(\frac{1-e^{-\varkappa t}}{\varkappa}\right)^2\sup_{[0,t]}|
    \widehat F_0(\tau)|
  \end{equation}

This  leads to
  \begin{equation}
    \label{eq:zero:1}
    |\widehat u_0(t)|^2\leq 2\left\{ \widehat f_{ 0}^2+ 
      g^2_{0}\left(\frac{1-e^{-2\varkappa 
            t}}{2\varkappa}\right)^2+\frac{1}{4}\left(\frac{1-e^{-\varkappa 
            t}}{\varkappa}\right)^4\left(\sup_{[0,t]}|
        \widehat F_0(\tau)|\right)^2\right\}.
  \end{equation}

  Summing up the inequalities \eqref{eq:e-1}, \eqref{eq:e-2},
  \eqref{eq:e-3} and \eqref{eq:zero:1} imply that inequality
  \eqref{eq:estimate-linear} holds, and this completes the proof of
  Proposition \ref{prop:3.2}.

\end{proof}

We turn now to prove the main result of this section, namely,  the proof of 
Theorem \ref{thm:1}.

\begin{proof}[Proof of Theorem \ref{thm:1} by a fixed point argument]

  Based on the energy estimate \eqref{eq:estimate-linear} of the
  solution to the linear Cauchy problem
  \eqref{eq:wave:7}--\eqref{eq:wave:8}, we shall show the existence of
  classical solutions to the Cauchy problem
  \eqref{eq:wave:1}--\eqref{eq:wave:2} in the interval $[0,\infty)$ in the
  Sobolev space $ H^{m+1}$ for $m> \frac{5}{2}$, under the assumption
  that the initial data, as well as $a(t,x)$, are sufficiently small. 
  In order to achieve this, we define a linear operator
  \begin{equation}
    \label{eq:linear-operator}
    \mathscr{L}: C([0,\infty); H^{m+1})\to  C([0,\infty); H^{m+1}), 
  \end{equation}
  as follows. 
  Let $u=\mathscr L (v)$ be the solution to the linear equation
  \begin{subequations}
    \begin{align}
      \label{eq:wave:12}
      &\partial_t^2 u+2\varkappa \partial_t u -\Delta u=\mathord{e^{-\varkappa t}} 
        a(t,x)(1+v)^3\\ 
      \label{eq:wave:13}
      & (u(0,x),\partial_t u(0,x))=(f(x),g(x)).
    \end{align}
  \end{subequations}
  Next, for $R>0$ we define a bounded set $B_R\subset H^{m+1}$ as
  follows
  \begin{equation}
    \label{eq:ball}
    B_R =\{v(t,\cdot)\in  C([0,\infty); H^{m+1}): 
    \sup_{[0,\infty)}\|v(t,\cdot)\|_{H^{m+1}}\leq R,\  
    v(0,x)=f(x), 
    \partial_t 
    v(0,x)=g(x)\}.
  \end{equation}

  Obviously, the ball $B_R$ is a closed set in the Banach space
  $C([0,\infty); H^{m+1})$, and that is why we can apply the Banach fixed
  point theorem to the operator $\mathscr{L}$, which will enable us to
  prove the existence of global solutions.
  In order to apply the Banach fixed point theorem we need to show:
  \begin{enumerate}[label=\alph*)]
    \item \label{item:section3-fourier:2}  $\mathscr L:B_R\to
    B_R$, that is, $\mathscr L$ maps the  ball into itself.
    \item \label{item:section3-fourier:3} $\mathscr L:B_R\to B_R$ is a
    contraction.
  \end{enumerate}

  We start with \ref{item:section3-fourier:2}:
 
 We shall use the energy estimate provided by Proposition
    \ref{prop:1}.
    So we set
    \begin{align*}
      M_1 & =\max\{2e^{-2\varkappa t}\left((1+2\varkappa^2(1+t^2)\right):t\geq 0\};\\ 
      M_2 & =\max\{4e^{-2\varkappa t}\left(1+t^2\right):t\geq 0\};\\
      M_3 & =\max\{e^{-2\varkappa t}\left(t^2(1+t^2\right)):t\geq 0\}.
    \end{align*}
    Using standard calculus in $H^m(\setT^3)$ there is a constant
    $C(R) $ such that
    \begin{equation}
      \|(1+v)^3\|_{L^\infty}\leq C_e \|(1+v)^3\|_{H^m}\leq C(R)
    \end{equation}
    for any $ v\in B_R$, here $C_e$ is the constant of the embedding
    $L^\infty \hookrightarrow H^m$.
    We first estimate the integral term and $\widehat F_0$ of the right
    hand side of \eqref{eq:estimate-linear}. 
    So
    \begin{equation}
      \int_0^t\|F(\tau,\cdot)\|_{\mathbullet H^m}^2d\tau\leq 
      t\sup_{[0,t]}\|a(\tau,\cdot)(1+v(\tau,\cdot)^3\|_{\mathbullet H^m}^2\leq t C_m^2
      \sup_{[0,\infty)}\|a(\tau,\cdot)\|_{\mathbullet H^m}^2C^2(R),
    \end{equation}
    where $C_m$ is the constant of the multiplication in the Sobolev
    space $H^m$. 
    Now,
    \begin{equation}
    \label{eq:est-0}
      \widehat F_0(\tau)=\frac{1}{(2\pi)^3}\int_{\setT^3}a(\tau,x)(1+v(\tau,x))^3dx,
    \end{equation}

    By Jensen's inequality (see e.~g.~\cite[Ch.~2]{Lieb-Loss}),
    \begin{equation}
      \begin{split}
        | \widehat F_0(\tau)|^2
        &=\frac{1}{(2\pi)^3}\int_{\setT^3}a^2(\tau,x)(1+v(\tau, x))^6dx\leq
        \|(1+v(\tau, \cdot))^3\|_{L^\infty}^2\|a(\tau,\cdot)\|_{L^2}^2 \\ & \leq
        C^2(R) \|a(\tau,\cdot)\|_{H^m}^2.
      \end{split}
    \end{equation}
    Now, by Proposition \ref{prop:3.2}, inequality
    \eqref{eq:estimate-linear},  $u=\mathscr L (v)$ satisfies the inequality
    \begin{equation}
      \begin{split}
        \|u(t,\cdot)\|_{H^{m+1}}^2 & \leq M_1\|f\|_{\mathbullet H^{m+1}}^2+
        M_2\|g\|_{\mathbullet H^{m}}^2+M_3C_m^2 C^2(R) \|a(\tau,\cdot)\|_{H^m}^2 \\
        & \widehat f_0^2+ \frac{\widehat g_0^2}{4\varkappa^2}+\frac{1}{4\varkappa^2} C^2(R)
        \|a(\tau,\cdot)\|_{H^m}^2 \leq R^2,
      \end{split}
    \end{equation}
    if
    \begin{equation}
      \label{eq:cond-1}
      \|f\|_{H^{m+1}}^2\leq \frac{ R^2}{4\max\{M_1,1\}},
    \end{equation}
    \begin{equation}
      \label{eq:cond-2}
      \|g\|_{H^{m}}^2\leq \frac{ R^2}{4\max\{M_2,\frac{1}{4\varkappa^2}\}}
    \end{equation}
    and
    \begin{equation}
      \label{eq:cond-3}
      \sup_{[0,\infty)}\|a(t,\cdot)\|_{H^{m}}^2\leq  
      \frac{R^2}{2C_m^2C^2(R)}\frac{1}{4\max\{M_3,\frac{1}{4\varkappa^4}\}}.
    \end{equation}
    Thus, $\mathscr L$ maps the ball into itself provided that
    \eqref{eq:cond-1}, \eqref{eq:cond-2} and \eqref{eq:cond-3} hold.

    \ref{item:section3-fourier:3} {Contraction:} Let
    $w=\mathscr L(v_1)-\mathscr L(v_2)$, then $w$ satisfies
    \begin{align}
      \label{eq:wave:14}
      &\partial_t^2 w+2\varkappa \partial_t w -\Delta u=\mathord{e^{-\varkappa t}} 
        a(t,x)\left((1+v_1)^3-(1+v_2)^3\right),\\ 
      \label{eq:wave:15}
      & (w(0,x),\partial_t w(0,x))=(0,0).
    \end{align}

    By the energy estimate \eqref{eq:estimate-linear}, we obtain that
    \begin{equation}
      \begin{split}
        \|w\|_{H^{m+1}}^2 & \leq 2 e^{-2\varkappa
          t}t^2(1+t^2)\sup_{[0,t]}\|a(\tau,\cdot)\left((1+v_1(\tau,
          \cdot))^3-(1+v_2(\tau, \cdot))^3\right)\|_{H^m}^2 \\ &+
        \frac{1}{4\varkappa^4}\sup_{[0,t]} |\widehat F_0(\tau)|^2.
      \end{split}
    \end{equation}
    Note that
    \begin{equation*}
      \left((1+v_1)^3-(1+v_2)^3\right)=(v_1-v_2)\left(3+3(v_1+v_2)+ 
        (v_1^2+v_1v_2+v_2^2)\right),
    \end{equation*}
    and that similar to equation \eqref{eq:est-0} we obtain
    \begin{equation}
      |\widehat F_0(\tau)|^2\leq 
      \|a(\tau,\cdot)\|_{H^m}^2\|(v_1-v_2)(\tau,\cdot)\left(3+3(v_1+v_2)+ 
        (v_1^2+v_1v_2+v_2^2)\right)(\tau,\cdot)\|_{L^\infty}^2.
    \end{equation}
    So by the embedding $L^\infty\hookrightarrow H^m$, the
    multiplication property of $H^m$ and the fact that
    $v_1,v_2\in B_R$, there exists a constant $K(R)$ such that
    \begin{equation}
      \|w(t,\cdot)\|_{H^{m+1}}^2\leq \max\{M_3,\frac{1}{4\varkappa^2}\}K^2(R)
      \sup_{[0,\infty)}\|(v_1-v_2)(t,\cdot)\|_{H^m}^2\sup_{[0,\infty)}\|a(t,
      \cdot)\|_{H^m}^2
    \end{equation}
    holds.
    Thus, the operator $\mathscr L:B_R\to B_R$ is a contraction provided
    that
    \begin{equation}
      \label{eq:cond-4}
      \sup_{[0,\infty)}\|a(t,\cdot)\|_{H^m}^2\leq \frac{1}{2}\frac{1}{ 
        \max\{M_3,\frac{1}{4\varkappa^2}\}K^2(R)}.
    \end{equation}
    So let $\epsilon$ be the minimum of the upper-bounds
    \eqref{eq:cond-1}--\eqref{eq:cond-4}, then the existence of a
    unique global solution follows from the application of the Banach
    fixed point theorem. 
    The solution belongs to the ball $B_R$, and therefore
    inequality \eqref{eq:limit-zero} holds. 
    That completes the proof of theorem \ref{thm:1}.

\end{proof}

\biblio

%% file: section4-sh.tex
\section{The wave equation as a modified symmetric hyperbolic system}
\label{symm-hyper-const}

In this section, we investigate the questions of global existence and
asymptotic decay of classical solutions to the Cauchy problem
\eqref{eq:wave:1}--\eqref{eq:wave:2} by using the theory of symmetric
hyperbolic systems and the corresponding energy estimates. 
Since the relativistic Euler equations can be written as a symmetric
hyperbolic system
(see~\cite{BK8
}), it will enable us, in the future, to couple the semi-linear
equation \eqref{eq:wave:1} to the Euler equations
\eqref{eq:Euler-Nordstrom:6}.

It is a well-known fact that wave equations can be cast into
symmetric hyperbolic form.
It turns out, however, that we need a modification of this standard
procedure, which we will outline in the next subsection.
With this new system at hand, we are able to prove results similar to
those in Section \ref{sec:semi-linear-wave}.
There are, however, some important differences between the results in
both sections which we have to point out.
We do not require that the initial data have to
be small, and we can, even, drop the term $e^{-\varkappa t}$ and yet
obtain  global existence.

However, since we rely, to a certain extend, on properties of the
$\mathbullet H^{m}$ spaces and since the right hand side of the wave
equation \eqref{eq:section2-field:1}, the term
$\left( 1+u \right)^{3}\not\in \mathbullet H^{m}$, we will perform a
projection on the wave equation that allows us to obtain a system of
an ordinary differential equation and a modified wave equation with a
right hand side, that does belong to $\mathbullet H^{m}$.

\subsection{The projection of the wave equation}
\label{sec:proj-wave-equat}

Based on our observations made in section \ref{sec:non-homog-sobol}
about norms of the spaces $H^{m}$ and $\mathbullet H^m$, and in particular the
orthogonal decomposition of $H^m$ \eqref{eq:norm:3}, we define the
orthogonal projection $P_0 : H^m \to \setR$, by
\begin{equation}
  \label{eq:section4-sh6}
  \widehat{u}_0=P_0(u).
\end{equation}
We denote by $u_h$ the complementary projection, that is,
\begin{equation}
  \label{eq:section4-sh7}
  u_h=\left( Id-P_0 \right)u.
\end{equation}
Since $u_h$ belongs to $\mathbullet H^m$, its norm is given by the formula
\begin{equation}
  \label{eq:section4-sh8}
  \left\Vert u_h \right\Vert_{H^m}^2 =\sum\limits_{k\neq0}^{} \left\vert k
  \right\vert^{2m} \left\vert \widehat{u}_{k} \right\vert
\end{equation}
and obviously
\(  \left\langle \widehat{u}_0,u_{h} \right\rangle_m=0\)
holds, where the inner product is given by equation \eqref{eq:norm:4}.

We apply now the projections $P_0$ and $\Id-P_0$ to the wave equation
\eqref{eq:section2-field:1}, that is, 
\begin{align}
  \label{eq:section4-sh10}
  P_0 \left\{  \partial_t^2 u +2\varkappa \partial_t u -\Delta u \right\}    &= 
P_0 \left\{ \mathord{e^{-\varkappa t}} a(t,x)(1+u)^3 \right\}
\end{align}
and
\begin{equation}
 \label{eq:section4-sh11}
  \left( \Id-P_0 \right)\left\{  \partial_t^2 u +2\varkappa \partial_t u -\Delta 
u
                                  \right\}=  \left( \Id -P_0 \right) \left\{ 
\mathord{e^{-\varkappa t}} a(t,x)(1+u)^3 \right\}.
\end{equation} 

Those projections result in the following system
\begin{subequations}
  \begin{align}
    \label{eq:section4-sh12}
    & \widehat{u}_0^{\prime\prime} + 2\varkappa \widehat{u}_0^{\prime} = 
      e^{-\varkappa  t}\widehat F_0 \\
    \label{eq:section4-sh32}
    &  \partial_t^2 u_h +2\varkappa \partial_t u_h -\Delta u_h     = 
      \mathord{e^{-\varkappa t}} \left( a(t,x)(1+u)^3 -\widehat{F}_{0} \right),
  \end{align}
\end{subequations}
where
\begin{equation}
  \widehat F_0= P_0(a(t,x) \left( 1+u 
  \right)^3)= \int_{\setT^3} a(t,x) \left( 1+u \right)^3dx.
\end{equation}

\subsection{A semi-linear wave equation written as a symmetric
  hyperbolic system} 
\label{sec:semi-linear-wave-1}

The most common way to write the wave equation as a symmetric hyperbolic
system is to consider either the vector valued function
\begin{equation*}
  V=
  \begin{pmatrix}
    \partial_t u\\
    \partial_x u
  \end{pmatrix}
  \qquad\mbox{or} \qquad
V=
  \begin{pmatrix}
    \partial_t u\\
    \partial_x u\\
    u
  \end{pmatrix}
\end{equation*}
as an unknown  (here $\partial_x u\eqdef (\partial_1 u,\partial_2 u,\partial_3 u)^{\intercal}$).
However, in both cases, for a system with  damping terms, the energy
estimates obtained are not appropriate  to show global existence.

We, therefore, introduce a different unknown by  setting
\begin{equation}
\label{eq:section5-sh:5}
  V\eqdef
  \begin{pmatrix}
    \partial_tu_h +\varkappa u_h \\
    \partial_x u_h
  \end{pmatrix}.
\end{equation} 
Then equation \eqref{eq:section4-sh32} can be written as a symmetric
hyperbolic system as follows
\begin{equation}
  \label{eq:symm:1}
  \begin{split}
    \partial_t V=\sum_{k=1}^3 B^k\partial_k V-\varkappa V+\varkappa^2
    \begin{pmatrix}
      u_{h} \\
      0 \\
      0\\
      0
    \end{pmatrix}
    +\mathord{e^{-\varkappa t}}
    \begin{pmatrix}
      a(t,x) (1+u)^3
      \\
      0 \\
      0 \\
      0
    \end{pmatrix}
    - e^{-\varkappa t}
    \begin{pmatrix}
      \widehat{F_{0}}
      \\
      0 \\
      0 \\
      0
    \end{pmatrix}
  \end{split}
\end{equation}
where $B^k$ are constant symmetric matrices,
\begin{equation*}
  B^k=\begin{pmatrix}
    0 & \delta_1^k & \delta_2^k & \delta_3^k \\
    \delta_1^k & 0 & 0 & 0\\
    \delta_2^k & 0 & 0 & 0\\
    \delta_3^k & 0 & 0 & 0
  \end{pmatrix}.
\end{equation*}

\subsection{Energy estimates}
\label{sec:energy-estimate}

\begin{defn}[The energy functional]
The energy functional for the unknown $V$, given by equation \eqref{eq:section5-sh:5} is
\begin{equation}
  \label{eq:energy:1}
  E(t)\eqdef \langle V(t),V(t)\rangle_m=\|\partial_t   u_h+\varkappa 
u_h\|_{H^m}^2+\|\partial_x u_h\|_{ H^m}^2.
\end{equation}
\end{defn}
\begin{rem}[About the definition of the energy]
\label{rem:section4-sh:1}
  It might look surprising to define the energy as a scalar product in
  $H^{m}$ while the vector $V$ as defined in \eqref{eq:section4-sh32}
  only contains terms that belong to $\mathbullet H^{m}$. 
  We use this notation since the energy estimates contain terms that
  belong to $H^m$.
\end{rem}

Proceeding in the usual way, suppose $V$ satisfies \eqref{eq:symm:1}, then  
differentiation of  the energy with respect to  
time results  that
\begin{equation*}
  \begin{split}
    \frac{1}{2}\frac{d}{dt} E(t) & = \langle \partial_t V(t),V(t)\rangle_m =\sum_{k=1}^3 \langle
    B^k\partial_k V, V\rangle_m -\varkappa \langle V,V
    \rangle_m+\varkappa^2 \langle u_h, \partial_t u_h+\varkappa u_h\rangle_m\\
    & + \mathord{e^{-\varkappa t}} \langle a(t,\cdot)(1+u)^3,\partial_t
    u_h+\varkappa u_h\rangle_m - \mathord{e^{-\varkappa t}} \langle
    \widehat F_0,\partial_t u_h+\varkappa u_h \rangle_m\\
    & = -\varkappa \| V\|_{H^m}^2+ \varkappa^2 \langle u_h, \partial_t u_h+\varkappa u_h\rangle_m
    + \mathord{e^{-\varkappa t}} \langle a(t,\cdot)(1+u)^3,\partial_t
    u_h+\varkappa u_h\rangle_m,
  \end{split}
\end{equation*}
 since $\mathbullet H^{m}\perp\setR $, and 
since $B^k$ are symmetric and constant, then by integration by
parts that $\langle B^k\partial_k V, V\rangle_m =0$.
By the Cauchy Schwarz inequality, we obtain
\begin{equation*}
  |\langle u_h, \partial_t u_h+\varkappa u_h\rangle_m|\leq   \|  u_h\|_{ H^m} 
  \| \partial_t u_h+\varkappa 
  u_h\|_{ H^m}
  \leq \| u_h\|_{ H^m}\| V\|_{ H^m}
\end{equation*}
and 
\begin{equation*}
  |\langle  a(1+u)^3,\partial_t    
u_h+\varkappa u_h\rangle_m|
    \leq \Vert a \left( 1+u \right)^{3}\Vert_{H^{m}}\| V\|_{ H^m},
\end{equation*}

which allows us to conclude, using the definition of the energy 
\eqref{eq:energy:1}, that 
\begin{equation}
  \label{eq:energy:2}
  \frac{1}{2}  \frac{d}{dt} E(t)\leq -\varkappa E(t)
  + \left\{ \varkappa^2 \left\Vert u_h(t) \right\Vert_{H^m} + e^{-\varkappa 
      t}\left\Vert a(t,\cdot)( \left( 1+u(t) \right)^{3} \right\Vert_{H^{m}} \right\}
  \sqrt{E(t)}
\end{equation}
We now apply Gronwall's inequality, Lemma \ref{lem:Gronwall}, in the
interval $[t_0,t]$ and with  $A(t)=-\varkappa$, then we obtain
\begin{equation}
  \label{eq:Gronwall:2}
  \begin{split}
    \sqrt{E(t)} & \leq e^{-\varkappa(t-t_0)}\sqrt{E(t_0)} +
    \varkappa^2\int_{t_0}^te^{-\varkappa(t-s)} \left\Vert u_h(s) \right\Vert_{H^{m}}ds
    \\ & + \int_{t_0}^te^{-\varkappa (t-s)}e^{-\varkappa s} \left\Vert
      a(s,\cdot) \left( 1+u(s) \right)^3 \right\Vert_{H^m}ds
  \end{split}.
\end{equation}

\begin{rem}[The role of $u_{h}$ in the a-priori estimates]
  \label{rem:section4-sh:2}
  We observe that the term $1+u=1+\widehat u_0+u_h$ implies that, for a
  fixed $\widehat u_0$ equation \eqref{eq:section4-sh32} is not coupled
  to \eqref{eq:section4-sh12} and it consists only of the unknown
  $u_h$. 
  Obviously, $u_h$ is not a solution to the system
  \eqref{eq:section4-sh12}--\eqref{eq:section4-sh32}, but it enables us
  to obtain important \textit{a-priori} estimates for the solution.
\end{rem}

We are now in a position to apply this energy estimate to show global
existence by a bootstrap argument, which is done in the next section.

\subsection{Global existence by a bootstrap argument}
\label{sec:glob-exist-bootstr}

In this section, we take the initial data in the homogeneous space,
\begin{equation}
 \label{eq:section5-sh:33}
\left\{\begin{array}{ll}
         u(x,0)=f(x), & \partial_t u(x,0)=g(x),\\
         f\in \mathbullet H^{m+1}, & g\in \mathbullet H^m.
       \end{array}\right.
\end{equation} 

This is not essential for the proof, but it makes it somewhat simpler. 
The following theorem is the main result of this section.

\begin{thm}[Global existence and decay of solutions]
  \label{thr:section5-sh:1}
  Let $0<\varkappa<1$ and $m> \frac{5}{2}$, let the initial data be as
  specified by \eqref{eq:section5-sh:33}, and $a\in C([0,\infty);H^{m})$.
  There exists a suitable constant $\varepsilon$ such that if
  \begin{equation}
    \label{eq:small:1}
    \sup_{[0,\infty)}\|a(t,\cdot)\|_{ H^{m}}<\varepsilon, 
  \end{equation}
  then system \eqref{eq:section4-sh12}--\eqref{eq:section4-sh32}, or
  equivalently equation \eqref{eq:wave:1}, with initial data given by
  \eqref{eq:section5-sh:33} has a unique solution
  \begin{equation}
    \label{eq:section4-sh:10}
    u\in C([0,\infty); H^{m+1}). 
  \end{equation}
  Moreover
  \begin{equation}
    \label{eq:decay}
    \lim_{t\to\infty}\|u(t)\|_{H^{m+1}}\leq \widetilde \epsilon,
  \end{equation}
  where $\widetilde \epsilon$ depends on the smallness condition
  \eqref{eq:small:1}.

\end{thm}

\begin{rem}[Comparison with theorem \ref{thm:1}]
  \label{rem:section4-sh:3}
  We emphasize that, contrary to Theorem \ref{thm:1}, in section
  \ref{subsec:homogeneous}, the initial data are not required to be small.
\end{rem}

\begin{rem}[About the asymptotic behavior of the metric]
\label{rem:section4-sh:5}
Recall that the physical metric has the following form
$ g_{\alpha\beta}=\phi^2\eta_{\alpha\beta}$, where
$\eta_{\alpha\beta} $ denotes the Minkowski metric. 
In Section \ref{sec:fields-equations} we concluded that the background
metric has the following form $(e^{\varkappa t})^2 \eta_{\alpha\beta}$.
We shall now use the asymptotic estimate \eqref{eq:decay} of the
global solutions to compare the asymptotic of the physical metric with
the background metric.  
We remind that $\phi=\mathord{e^{\varkappa t}}(1+u)$, where $u$ is the
solution to the Cauchy problem \eqref{eq:wave:1}--\eqref{eq:wave:2},
therefore we can conclude that the asymptotic behavior of the metric
can be described by the following expression
  \begin{equation}
    ( 1-\widetilde \epsilon)^2\leq 
    \lim_{t\to\infty}\frac{g_{\alpha\beta}(t,x)}{e^{2\varkappa t}\eta_{\alpha\beta} 
    }=
    \lim_{t\to\infty}(1+u(t,x))^2\leq (1+\widetilde \epsilon)^2.
  \end{equation}
   
\end{rem}

The proof of Theorem \ref{thr:section5-sh:1} is based on the following
propositions which we present together with their corresponding
proofs.
We recall that by the existence theorem, Theorem
\ref{thr:section3-fourier:1}, the solution of the system
\eqref{eq:section4-sh12}--\eqref{eq:section4-sh32} exists in a certain
time interval $[0,T]$.

\begin{prop}[A priori estimates]
  \label{prop:2}
  Let $ 0<\varkappa<1$, $1<\beta<\frac{1}{\varkappa}$ and set $\alpha= \|V(0)\|_{H^m}$.
  Assume the  solution $u=\widehat u_0+u_h$ to 
  \eqref{eq:section4-sh12}--\eqref{eq:section4-sh32} with initial data
  \eqref{eq:section5-sh:33} exists for $t\in [0,T]$. 
  If $\mathord{\left\|{a(t,\cdot)}\right\|_{ H^m}}$ is sufficiently small,
  then there exists a $T^{+}$, $0<T^{+}\leq T$, such that
  \begin{equation}
    \label{eq:section4-sh:2}
    \sup_{[0,T^{+}]}\mathord{\left\|{u(t)}\right\|_{ H^m}}\leq \alpha\beta
  \end{equation}
  and
  \begin{equation}
    \label{eq:section4-sh:1}
    E(T^{+})\leq E(0).
  \end{equation}
\end{prop}

\begin{proof}
  We start with the proof of inequality \eqref{eq:section4-sh:2}. 
  Recall that although we have written $\alpha= \|V(0)\|_{H^m}$, the initial
  data are in the homogeneous Sobolev space and therefore by
  Proposition \ref{prop:9} and \eqref{eq:energy:1}, we obtain
  \begin{equation}
    \|u(0)\|_{H^{m+1}}=\|f\|_{H^{m+1}}=\|\partial_x f\|_{H^{m}}\leq 
    \|V(0)\|_{H^m}=\alpha.
  \end{equation}
  Hence, since $ \beta > 0$, it follows from the existence Theorem
  \ref{thr:section3-fourier:1} and the continuity property of the
  corresponding solutions, that there exists $0<T^+\leq T$ such that
  \begin{equation}
    \label{eq:section4-sh4}
    \sup_{[0,T^+]}\| u(t)\|_{ H^m}\leq \alpha\beta.
  \end{equation}
  We now turn to inequality \eqref{eq:section4-sh:1}.
  For $t\in [0,T^+]$ we observe, using inequality \eqref{eq:Gronwall:2}
  that
  \begin{equation}
    \label{eq:Gronwall:1}
    \begin{split}
      \sqrt{E(t)} & \leq e^{-\varkappa t}\sqrt{E(0)} +
      \varkappa^2\int_{0}^te^{-\varkappa (t-s)} \alpha\beta ds +
      \int_{0}^te^{-\varkappa (t-s)}e^{-\varkappa s} \left\Vert a(s,\cdot)
        \left( 1+u(s) \right)^3 \right\Vert_{H^m} ds\\ & \leq e^{-\varkappa
        t}\sqrt{E(0)} +\varkappa \left(1 -e^{-\varkappa t}\right)\alpha\beta+
      te^{-\varkappa t}\sup_{[0,t]}\left\Vert a(s,\cdot) \left( 1+u(s)
        \right)^3 \right\Vert_{H^m}.
    \end{split}
  \end{equation}
  A simple algebraic manipulation shows us that, $E(t)\leq E(0)$, if
  \begin{equation}
    \varkappa \left(e^{\varkappa t}-1\right)\alpha\beta+  
    t\sup_{[0,t]}\left\Vert 
      a(s,\cdot) \left(1+u(s) \right)^3 \right\Vert_{H^m}\leq 
    \left(e^{\varkappa 
        t}-1\right)\sqrt{E(0)},
  \end{equation}
  or equivalently
  \begin{equation}
    \label{eq:section4-sh:7}
    t\sup_{[0,t]}\left\Vert 
      a(s,\cdot) \left(1+u(s) \right)^3 \right\Vert_{H^m}\leq 
    \left(e^{\varkappa 
        t}-1\right)\alpha\left(\beta-\varkappa\right).
  \end{equation}
  Since for $s\in [0,T^+]$, we can conclude that
  $\|u(s)\|_{H^m}\leq \alpha\beta\leq 2\alpha \beta$, we can apply
  Proposition \ref{prop:1} with $A=2\alpha \beta$, that results in
  \begin{equation}
    \label{eq:energy:3}
    \left\Vert 
      a(s,\cdot) \left(1+u(s) \right)^3 \right\Vert_{H^m}\leq 
    C(2\alpha\beta)\|a(s,\cdot)\|_{H^m}.
  \end{equation}
  We now set
  \begin{equation}
    \label{eq:energy:4}
    \epsilon_0=\frac{\varkappa\alpha(\beta-\varkappa)}{C(2\alpha\beta)}.
  \end{equation}
  Since $\beta-\varkappa>0$, $\epsilon_0>0$, therefore we can demand the smallness
  condition
  \begin{equation}
    \label{eq:energy:5}
    \sup_{[0,\infty)}\|a(t,\cdot)\|_{H^m}\leq \epsilon_0.
  \end{equation}
  We now let $t=T^+$ in inequality \eqref{eq:section4-sh:7}, then by
  inequality \eqref{eq:energy:3}, condition \eqref{eq:energy:5}, with
  \eqref{eq:energy:4}, we conclude that
  \begin{equation}
    \label{eq:section4-sh3}
    \begin{split}
      T^+\sup_{[0,T^+]}\left\Vert a(s,\cdot) \left(1+u(s) \right)^3 \right\Vert_{H^m} &
      \leq T^+ C(2\alpha\beta)\sup_{[0,T^+]}\|a(s,\cdot)\|_{H^m}\leq T^+
      C(2\alpha\beta)\epsilon_0 \\ & = \varkappa T^+\alpha(\beta-\varkappa)\leq (e^{\varkappa
        T^+}-1)\alpha(\beta-\varkappa),
    \end{split}
  \end{equation}
  holds and consequently \eqref{eq:section4-sh3} implies inequality 
  \eqref{eq:section4-sh:7}.
  This proves \eqref {eq:section4-sh:1} and completes the proof of the
  proposition. 
  In the last step, we used the elementary inequality $x\leq e^{x}-1$.
\end{proof}

Based on  Proposition \ref{prop:2} we define 
\begin{defn}[Definition of $T^\star$]
\begin{equation}
    \label{eq:section4-sh:4}
    T^\ast=\sup\left\{T: \sup_{[0,T]}\mathord{\left\|{u(t)}\right\|_{H^m}}\leq 
  \alpha\beta  \ \text{and} \ E(T)\leq E(0)\right\}.
  \end{equation}
\end{defn}
The following proposition plays a central role in proving Theorem 
\ref{thr:section5-sh:1}. 
\begin{prop}[$T^{\star}$ is not finite]
  \label{prop:10}
 Under the assumptions  of Proposition \ref{prop:2}, we obtain
  \begin{equation*}
    T^\ast=\infty.
  \end{equation*}
\end{prop}

It is important to note that we need two conditions in the definition
of $T^{\star}$, as Proposition \ref{prop:2} already suggests. 
The role of these two conditions will become clearer after we finish 
the proof and we will come back to this point.

\textsc{Sketch of the proof:} The proof of this proposition is rather
long, and as we said, crucial for theorem \ref{thr:section5-sh:1} and
that is why we sketch here its structure.
We prove Proposition \ref{prop:10} by a contradiction argument, in
other words, we assume that $T^\ast$ is finite, and then we show that
both conditions of \eqref{eq:section4-sh:4} hold in a larger interval.

The first step of the proof deals with the extension of the solution for $t>T^\ast$.
In the second step, using the inequality \eqref{eq:section4-sh:1}, we show that 
there exists  a $T^{\ddagger} >T^\ast$ such that \begin{equation}
  \label{eq:section4-sh:3}
  \sup_{[0,T^\ddagger]}\mathord{\left\|{u(t)}\right\|_{H^m}}\leq\alpha\beta.  
\end{equation}
With this inequality proven, we are then able, in the third step, to
show that
\begin{equation}
  \label{eq:section4-sh14}
  E(T^\ddagger)\leq E(0).    
\end{equation}

The existence of these inequalities in the interval $[0,T^{\ddagger}]$
contradicts the definition of $T^{\ast}$. 
Therefore we conclude that $T^\ast=\infty$.

\begin{proof}[Proof of Proposition \ref{prop:10}] \quad
  \begin{enumerate}[label=\textsc{Step \arabic*.},wide,labelwidth=!,labelindent=0pt]
    \item\label{item:section4-sh:3}
    We need to extend the solution beyond $T^{\star}$, so let
    $\widetilde u$ be the solution to equation \eqref{eq:wave:1} with
    initial data $\widetilde u(T^\ast,x)=\widetilde f(x)$ and
    $ \partial_t \widetilde u(T^\ast,x)=\widetilde g(x)$, where
    $\widetilde f(x) =u(T^\ast,x)$ and
    $\widetilde g(x)=\partial_t u(T^\ast,x)$.

    We will show that
    \begin{equation}
      \label{eq:outline-of-proof-prop5:1}
      u(t)\in H^{m+1} \quad \mbox{for}\quad t\in [0,T^\ast], 
    \end{equation}
    which implies $ \widetilde f\in H^{m+1}$, a fact that is needed in
    order to apply the existence theorem, Theorem
    \ref{thr:section3-fourier:1}.
    To prove \eqref{eq:outline-of-proof-prop5:1} we apply Proposition
    \ref{prop:9} to $ u_{h} $, which allows us to conclude
    \begin{equation}
      \label{eq:section4-sh35}
      \|u(t)\|_{H^{m+1}}^2=|\widehat u_0(t)|^2+\|u_h(t)\|_{H^{m+1}}^2=|\widehat u_0(t)|^2+\|\partial_x u_h(t)\|_{H^{m}}^2.
    \end{equation}
    Hence, since $u_h(t)\in H^m$, we see by equation
    \eqref{eq:section4-sh35} that $\partial_x u_h(t)\in H^{m}$ and we can conclude
    $u(t)\in H^{m+1}$. 
    Consequently, by the existence theorem, Theorem
    \ref{thr:section3-fourier:1}, there exists a $T_1>T^\ast$ such that
    $\widetilde u(t)$ exists for $t\in [T^\ast,T_1]$.

    \item\label{item:section4-sh:4}
    We turn now to the proof of inequality \eqref{eq:section4-sh:3}. 
    Since $\|\widetilde u(T^\ast)\|_{H^m}\leq \alpha\beta$, there exits a
    $T_{2}$, $T^\ast<T_2\leq T_1$, such that
    \begin{equation}
      \label{eq:outline-of-proof-prop54}
      \sup_{ [T^{\star},T_{2}]}\|\widetilde u(\tau)\|_{H^m}\leq 2\alpha\beta.
    \end{equation}
    holds. 
    We now set
    \begin{equation}
      \label{eq:tilde}
      \widetilde V=\begin{pmatrix} \partial_t \widetilde u_h+\varkappa
        \widetilde u_h\\
        \partial_x \widetilde u_h
      \end{pmatrix},
    \end{equation}
    then $(\widehat{\widetilde u}_0,\widetilde u_h)$ solves system
    \eqref{eq:section4-sh12}--\eqref{eq:section4-sh32} with the initial
    data
    \begin{subequations}
      \begin{align}
        & \widehat{\widetilde u}_0(T^\ast)=\widehat{\widetilde f}_0, \quad
          \partial_t\widehat{\widetilde u}_0(T^\ast)=\widehat{\widetilde g}_0
        \\
        & \widetilde u_h(T^\ast,x)=\widetilde f_h(x), \quad  
          \partial_t\widetilde u_h(T^\ast,x)=\widetilde g_h(x).
      \end{align}
    \end{subequations}

    We shall estimate each component of
    $(\widehat{\widetilde u}_0,\widetilde u_h)$ separately. 
    We take $0<\epsilon_1$ such that $1<\beta -\epsilon_1$ and then we will prove below
    the following two inequalities:
    \begin{equation}
      \label{eq:outline-of-proof-prop53}
      |\widehat{\widetilde u_0}(t)|\leq \frac{\alpha \beta\epsilon_1}{2}, \qquad       t\in [T^\ast,T_4],
    \end{equation}
    and
    \begin{equation}
      \label{eq:outline-of-proof-prop51}
      \sup_{[T^\ast,T_5]} \|\widetilde u_h(t)\|_{H^{m}}\leq \alpha(\beta-\epsilon_1). 
    \end{equation}
    Here $T_4,T_5\in (T^\ast,T_2]$. 
    Setting $T^{\ddagger}=\min\{T_4,T_5\}$, and combining
    \eqref{eq:outline-of-proof-prop53} and
    \eqref{eq:outline-of-proof-prop51}, we obtain
    \begin{equation}
      \label{eq:section4-sh:5}
      \sup_{[T^\ast,T^{\ddagger}]} \|\widetilde u(t)\|_{H^m}^2=  
      \sup_{[T^\ast,T^{\ddagger}]}|\widehat{\widetilde u_0}(t)|^2+  
      \sup_{[T^\ast,T^{\ddagger}]}\|\widetilde 
      u_h(t)\|_{H^m}^2\leq\frac{(\alpha\beta\epsilon_1)^2}{4}
      +(\alpha(\beta-\epsilon_1))^2\leq (\alpha\beta)^2,
    \end{equation}
    which proves \eqref{eq:section4-sh:3}. 
    The last inequality requires that $\epsilon_1 \leq \frac{8}{\beta+4}$.

    We start with the estimate of $\widehat{\widetilde u_0}(t)$ for
    $t\in [T^\ast,T_2]$. 
    Since it satisfies the initial value problem
    \eqref{eq:fs:4}--\eqref{eq:fs:5} with $k=0$, its solution is given
    by
    \begin{equation}
      \label{eq:section4-sh:6}
      \widehat{\widetilde u}_0(t) = \widehat{\widetilde 
        f}_0+\widehat{\widetilde 
        g}_0\left(\frac{1-e^{-2\varkappa(t-T^\ast)}}{2\varkappa}\right) +
      \frac{1}{2\varkappa}\int_{T^\ast}^t\left(1-e^{
          -2\varkappa(t-\tau) }
      \right)e^{-\varkappa \tau} \widehat F_{ 0}(\tau) d\tau,
    \end{equation}
    where
    \begin{equation}
      \widehat F_{ 
        0}(\tau)=\frac{1}{(2\pi)^3}\int_{\setT^3}a(\tau,x)(1+\widetilde 
      u(\tau,x))^3 dx.
    \end{equation}
    Since $\widehat{\widetilde f}_0=\widehat u_0(T^\ast)$ and the initial
    data $f,g\in \mathbullet H^{m}$, we conclude by equation  \eqref{eq:u-0} and
    estimate \eqref{eq:u-0:2} that
    \begin{equation}
      \label{eq:estimate:7}
      | \widehat{\widetilde f}_0|\leq \frac{1}{2}\left(\frac{1-e^{-\varkappa 
            T^{\ast}}}{\varkappa}\right)^2\sup_{[0,T^\ast]}|\widehat F_0(t)|
    \end{equation}
    holds. 
    By Proposition \ref{prop:1}, we obtain
    \begin{equation}
      \label{eq:estimate:8}
      |\widehat F_0(t)|\leq 
      \sup_{[0,T^\ast]}\|a(t,\cdot)\|_{L^\infty}C(2\alpha\beta).
    \end{equation}
    Hence, if we require that
    \begin{equation}
      \label{eq:cond:1}
      \sup_{[0,\infty)}\|a(t,\cdot)\|_{L^\infty}\leq \frac{2\alpha\beta\varkappa^2\epsilon_1}{6 C(2\alpha\beta)},
    \end{equation}
    then we conclude that
    \begin{equation}
      \label{eq:estimate:1}
      | \widehat{\widetilde f}_0|\leq \frac{\alpha\beta\epsilon_1}{6}.
    \end{equation}
    holds. 
    
    Before we proceed, we remark, that the smallness condition
    \eqref{eq:cond:1} is given in the term of $L^\infty$ norm, but it
    can easily be formulated in terms of $H^m$ norm by Sobolev
    embedding theorem.

    Next, since
    $\displaystyle\lim_{t\to
      T^\ast}\left(\frac{1-e^{-2\varkappa(t-T^\ast)}}{2\varkappa}\right)=0$,
    there exists a $T_3$, with $T^\ast<T_3\leq T_2$, such that
    \begin{equation}
      \label{eq:estimate:2}
      \left|\widehat{\widetilde g_0}\right| 
      \left(\frac{1-e^{-2\varkappa(t-T^\ast)}}{2\varkappa}\right)\leq 
      \frac{\alpha\beta\epsilon_1}{6}, \quad t\in [T^\ast,T_3].
    \end{equation}
    For the third term on the left side of equation
    \eqref{eq:section4-sh:6}, we use inequality
    \eqref{eq:outline-of-proof-prop54}, and then by a similar argument
    used to estimate the term $| \widehat{\widetilde f}_0|$, we obtain
    \begin{equation}
      \label{eq:estimate:6}
      \left| \frac{1}{2\varkappa}\int_{T^\ast}^t\left(1-e^{
            -2\varkappa(t-\tau) }
        \right)e^{-\varkappa \tau} \widehat F_{ 0}(\tau) d\tau\right|\leq     
      \frac{1}{2\varkappa^2}C(2\alpha\beta)\sup_{[0,\infty)}\|a(t,\cdot)\|_{L^\infty}
      \leq \frac{\alpha\beta \epsilon_1}{6},
    \end{equation}
    provided that the condition \eqref{eq:cond:1} is satisfied.
    Letting $T_4=\min\{T_2,T_3\}$, we conclude from inequalities
    \eqref{eq:estimate:1}, \eqref{eq:estimate:2} and
    \eqref{eq:estimate:6} implies that \eqref{eq:outline-of-proof-prop53}
    holds.

    We now turn to prove inequality
    \eqref{eq:outline-of-proof-prop51}: For a fixed
    $\widehat{\widetilde u_0}(t)$, the unknown $\widetilde V(t)$, as
    defined in equation \eqref{eq:tilde}, satisfies the symmetric hyperbolic
    system \eqref{eq:symm:1}, and that is why we can apply Proposition
    \ref{prop:6} with initial data $V(T^\ast,x)$. 
    We first observe by Proposition \ref{prop:9} that
    \begin{equation}
      \label{eq:outline-of-proof-prop55}
      \begin{split}
        \mathord{\left\|{\widetilde u_h(t)}\right\|_{H^m}} &\leq \|\partial_x
        \widetilde u_h(t)\|_{H^{m-1}}\leq \|\widetilde V(t)\|_{H^{m-1}} \leq \|
        \widetilde V(t) -\widetilde V(T^\ast)\|_{H^{m-1}}+ \| \widetilde
        V(T^\ast)\|_{H^{m-1}}
      \end{split}
    \end{equation}
    holds. 
    By the definition of $T^\ast$, we can conclude that
    $E(T^\ast)\leq E(0)$ holds, which then implies the inequality
    $\|V(T^\ast,\cdot)\|_{H^m} =\sqrt{ E(T^\ast)}\leq \sqrt{E(0)}=\alpha$. 
    We now can apply Proposition \ref{prop:6} to
    $\| \widetilde V(t) -\widetilde V(T^\ast)\|_{H^{m-1}}$ with $C_0$
    depending on $\alpha$, combine it with inequality
    \eqref{eq:outline-of-proof-prop55}, and then obtain
    \begin{equation}
      \label{eq:section4-sh:13}
      \mathord{\left\|{\widetilde u_h(t)}\right\|_{H^m}}  \leq 
      C_0(\alpha)(t-T^\ast)^{\frac{1}{m}}+ \alpha\leq
      \alpha(\beta-\epsilon_1)
    \end{equation}
   provided that
    $t-T^\ast\leq
    \left(\frac{\alpha(\beta-\epsilon_1-1)}{C_0(\alpha)}\right)^m$. 
    Thus \eqref{eq:outline-of-proof-prop51} holds with
    $T_5=T^\ast+\left(\frac{\alpha(\beta-\epsilon_1-1)}{C_0(\alpha)}\right)^m$.

    \item\label{item:section4-sh:5} It remains to show that
    $E(t)\leq E(0)$ for $t\in [T^\ast,T^{\ddagger}]$, where
    $T^\ddagger=\min\{T_4,T_5\}$.
    We will first establish the inequality
    \begin{equation}
      \label{eq:section4-sh:8}
      \sqrt{E(t)}            \leq e^{-\varkappa(t-T^\ast)}\sqrt{E(T^\ast)}+
      \varkappa\left(1-e^{-\varkappa(t-T^\ast)}\right)\alpha\beta+
      e^{-\varkappa t}(t-T^\ast)\sup_{[T^\ast,t]}\left\Vert
        a(\tau,\cdot) \left( 1+\widetilde u(\tau) \right)^3 \right\Vert_{H^m}.
    \end{equation}
    Using the energy estimate \eqref{eq:Gronwall:2} we observe that
    \begin{equation}
      \label{eq:Gronwall:8}
      \begin{split}
        \sqrt{E(t)} & \leq e^{-\varkappa(t-T^\ast)}\sqrt{E(T^\ast)} +
        \varkappa^2\int_{T^\ast}^te^{-\varkappa(t-\tau)}
        \left\Vert \widetilde u_h(\tau) \right\Vert_{H^{m}}d\tau \\
        & + \int_{T^\ast}^te^{-\varkappa (t-\tau)}e^{-\varkappa \tau}
        \left\Vert a(\tau,\cdot) \left( 1+\widetilde u(\tau) \right)^3 \right\Vert_{H^m}
        d\tau.
      \end{split}
    \end{equation}
    Since the inequality $\|\widetilde u_h(t)\|\leq \alpha\beta$ holds in the interval
    $[T^\ast,T^\ddagger]$, then inequality \eqref{eq:section4-sh:8} follows by
    inserting this bound into the energy estimate \eqref{eq:Gronwall:8}.
    From this inequality, we observe that
    \begin{math}
      \sqrt{E(t)}\leq \sqrt{E(0)}
    \end{math} holds if
    \begin{equation}
      \begin{split}
        & \varkappa\left(e^{\varkappa(t-T^\ast)}-1\right)\alpha\beta+
        e^{-\varkappa
          t}e^{\varkappa(t-T^\ast)}(t-T^\ast)\sup_{[T^\ast,t]}\left\Vert
          a(\tau,\cdot) \left(1+\widetilde u(\tau) \right)^3 \right\Vert_{H^m}\\ \leq &
        e^{\varkappa(t-T^\ast)}\sqrt{E(0)}-\sqrt{E(T^\ast)}=\sqrt{E(0)}\left(e^{
            \varkappa(t-T^\ast) }-1\right)+\sqrt{E(0)}-\sqrt{E(T^\ast)}.
      \end{split}
    \end{equation}
  
    But we already know that $E(T^\ast)\leq E(0)$ holds, by the
    definition of $T^\ast$. 
    Therefore it suffices to show that
    \begin{equation}
      \label{eq:outline-of-proof-prop512}
      \varkappa\left(e^{\varkappa(t-T^\ast)}-1\right)\alpha\beta+ 
      e^{-\varkappa T^\ast}(t-T^\ast)\sup_{[T^\ast,t]}\left\Vert 
        a(\tau,\cdot) \left(1+\widetilde u(\tau) \right)^3 \right\Vert_{H^m} 
      \leq 
      \sqrt{E(0)}\left(e^{\varkappa(t-T^\ast)}-1\right).
    \end{equation}
    Note, that since $\varkappa$ is strictly positive, we conclude inequality
    $e^{-\varkappa T^\ast}<1$, and therefore we can drop the term
    $e^{-\varkappa T^\ast}$. 
    So it is enough to show that
    \begin{equation}
      \label{eq:section4-sh:9}
      (t-T^\ast)\sup_{[T^\ast,t]}\left\Vert 
        a(\tau,\cdot) \left(1+\widetilde u(\tau) \right)^3 \right\Vert_{H^m} 
      \leq 
      \left(e^{\varkappa(t-T^\ast)}-1\right)\alpha(1-\beta\varkappa).
    \end{equation}
    We now let $t=T^{\ddagger}$ and we proceed as we did in the proof
    of Proposition \ref{prop:2}. 
    Under the smallness condition on $a(t)$, \eqref{eq:energy:5}
    where $\epsilon_0 $ is given by \eqref{eq:energy:4},  we conclude that
    \begin{equation}
      \label{eq:outline-of-proof-prop510}
      \begin{split}
        (T^{\ddagger}-T^\ast)\sup_{[T^\ast,T^{\ddagger}]}\left\Vert
          a(\tau,\cdot) \left(1+\widetilde u(\tau) \right)^3 \right\Vert_{H^m}
        &\leq(T^{\ddagger}-T^\ast)C(2\alpha\beta)\sup_{[0,\infty)}\|a(\tau,\cdot)\|_{H^m
        }\\
        \leq \varkappa (T^{\ddagger}-T^\ast)\alpha(1-\beta\varkappa)
        &\leq\left(e^{\varkappa(T^{\ddagger}-T^\ast)}-1\right)\alpha(1-\beta\varkappa).
      \end{split}
    \end{equation}
    holds. 
    We observe that the inequalities in expression
    \eqref{eq:outline-of-proof-prop510} imply inequality
    \eqref{eq:section4-sh:9}, and thus we have proved inequality
    \eqref{eq:section4-sh14}.         

    Taking into account the conditions on $a(t)$, namely
    \eqref{eq:energy:5} and \eqref{eq:cond:1} respectively, we
    conclude that there exists a positive $\epsilon$ depending on
    $\epsilon_0$ and $\epsilon_{1}$ such that if
    $\sup_{[0,\infty)}\|a(t,\cdot)\|_{H^m}\leq \epsilon$, the inequalities
    $\sup_{[0,T^{\ddagger}]}\left\{ \|u(t)\|_{H^m} \right\}\leq \alpha
    \beta$ and $E(T^\ddagger)\leq E(0)$ hold for $t\in [0, T^{\ddagger}]$. 
    It is important to note that the condition \eqref{eq:cond:1} also can
    be formulated in terms of the $H^m $ norm.
    Therefore, both conditions hold in the larger time interval
    $[0,T^{\ddagger}]$. 
    This implies that the assumption that $T^\ast<\infty$ is false and
    that completes the proof of Proposition \ref{prop:10}.
  \end{enumerate}
\end{proof}

\begin{rem}[About the definition of $T^\star$]
  We come back to the question of why we had two conditions in the
  definition of $T^{\star}$. 
  One motivation was Proposition \ref{prop:2}, however there is an
  important difference in the proofs of Proposition \ref{prop:2} and
  Proposition \ref{prop:10}. 
  While we proved
  \begin{equation}
    \sup_{[0,T]}\mathord{\left\|{u(t)}\right\|_{H^m}}\leq   \alpha\beta
  \end{equation}
  in Proposition \ref{prop:2} by a simple continuity argument, we
  needed condition $E(T^{\ddagger})\leq E(0)$ to prove that
  corresponding inequality in Proposition \ref{prop:10}.
  In other words, both conditions are interconnected appropriately.  
\end{rem}
\vspace{-0.35cm}
We turn now to the proof of Theorem \ref{thr:section5-sh:1}.
\vspace{-0.35cm}
\begin{proof}[Proof of Theorem \ref{thr:section5-sh:1}]
  The global existence and regularity essentially follows from
  Propositions \ref{prop:2} and \ref{prop:10}. 
  The asymptotic behavior \eqref{eq:decay} will be proven by
  considering $\left\vert \widehat{u}_{0} \right\vert$ and
  $\left\Vert u_{h} \right\Vert$ separately.
 
  We start to prove equation \eqref{eq:section4-sh:10}. 
  By the existence theorem, Theorem \ref{thr:section3-fourier:1}, the
  solution to the initial value problem
  \eqref{eq:wave:1}--\eqref{eq:wave:2} exists in a certain time
  interval $[0,T]$. 
  Consequently the system
  \eqref{eq:section4-sh12}--\eqref{eq:section4-sh32} has a solution in
  the time interval $[0,T]$. 
  Then we can apply Propositions \ref{prop:2} and \ref{prop:10} that
  provide the existence of a global solution $u\in C([0,\infty);H^m)$. 
  By Proposition \ref{prop:9} we obtain
  \begin{equation}
    \label{eq:outline-of-proof-theorem3:1}
    \|u(t)\|_{H^{m+1}}^2= |\widehat u_0(t)|^2+\|u_h(t)\|_{H^{m+1}}^2=
    |\widehat u_0(t)|^2+\|\partial_x u_h(t)\|_{H^{m}}^2.
  \end{equation}
  Hence, since $V\in C([0,\infty);H^m)$, it follows that $u\in C([0,\infty);H^{m+1})$.

  We now turn to the proof of the asymptotic behavior of the global
  solution as described by equation \eqref{eq:decay}. 
  The idea is again based on the decomposition $u=\widehat{u}_{0}+u_{h}$
  and to show that $\lim_{t\to\infty} \|u_h(t)\|_{H^{m+1}}=0$. So 
  we set
  \begin{equation}
    \label{eq:outline-of-proof-theorem3:4}
    \mu\eqdef\limsup_{t\to\infty}\|u_h(t)\|_{H^{m+1}}.
  \end{equation}
 Then    for a given $\epsilon >0$, there exists a $t_0$ such that
  $ \sup_{[t_0,\infty)} \|u_h(t)\|_{H^{m+1}}\leq \mu+\epsilon $.
  Using the energy estimate \eqref{eq:Gronwall:2} for $t>t_0$ and
  Proposition \ref{prop:1} with $A=\mu+\epsilon$, we obtain
  \begin{equation}
    \label{eq:Gronwall:13}
    \begin{split}
      \sqrt{E(t)}&\leq e^{-\varkappa(t-t_0)}\sqrt{E(t_0)}+
      \varkappa\left(1-e^{-\varkappa(t-t_0)}\right)\sup_{[t_0,t]} \|
      u_h(\tau)\|_{\mathbullet H^m} \\
      & \quad + \mathord{e^{-\varkappa t}}
      (t-t_0)\sup_{[t_0,t]}\left\{\|a(\tau,\cdot) (1+u(\tau))^3\|_{H^m}\right\}
      \\
      & \leq e^{-\varkappa(t-t_0)}\sqrt{E(t_0)}+
      \varkappa\left(1-e^{-\varkappa(t-t_0)}\right)(\mu+\epsilon)\\
      &\quad + \mathord{e^{-\varkappa t}}
      (t-t_0)C(\mu+\epsilon)\sup_{[0,\infty)}\|a(t,\cdot)\|_{H^m}.
    \end{split}
  \end{equation}
  We conclude from inequality \eqref{eq:Gronwall:13} that
  \begin{equation}
    \label{eq:asymp:1}
    \limsup_{t\to\infty}\sqrt{E(t)}\leq \varkappa(\mu+\epsilon)
  \end{equation}
 holds.
  On the other hand, using the fact that
  \begin{math}
    \|u_h(t)\|_{H^{m+1}}=\|\partial_x u_h(t)\|_{H^m}\leq \sqrt{E(t)}
  \end{math}
  we obtain
  \begin{equation}
    \label{sec:step-4}
    \mu\leq \limsup_{t\to\infty}{\sqrt{E(t)}}\leq \varkappa(\mu+\epsilon).
  \end{equation}
  We may assume that $\mu$ is strictly positive since otherwise there is nothing to be
  proven.  
  Since $\varkappa$ is strictly smaller than $1$, we can choose $\epsilon$ to be  $\epsilon=\left( 1-\varkappa \right)\mu>0$. 
  Then, from inequality \eqref{sec:step-4}, we obtain
  \begin{equation}
    \label{eq:section4-sh:15}
    \mu\leq \varkappa \left( \mu+\epsilon \right)=\varkappa\mu + \varkappa \left( 1-\varkappa \right)\mu<\varkappa\mu+ \left( 1-\varkappa
    \right)\mu=\mu,
  \end{equation}
  which implies that $\mu=0$ holds. 
  Recalling definition \eqref{eq:outline-of-proof-theorem3:4}, we
  conclude that
  \begin{equation}
    \label{eq:section4-sh:12}
    \lim_{t\to\infty}\|u(t)\|_{H^{m+1}}^2=\lim_{t\to\infty}|\widehat 
    u_0(t)|^2+\lim_{t\to\infty} \|u_h(t)\|_{H^{m+1}}^2=\lim_{t\to\infty}|\widehat 
    u_0(t)|^2,
  \end{equation}
  holds and that is why it remains to estimate the limit of the term
  $|\widehat u_0(t)|$ only.
  Since we use the following initial data, $f,g\in \mathbullet H^{m}$, we can
  express $\widehat{u}_0(t)$ explicitly by formula \eqref{eq:u-0} and
  observe that
  \begin{equation}
    \widehat u_{ 0}(t)=\frac{1}{2\varkappa}\int_0^t\left(1-e^{
        -2\varkappa(t-\tau) }
    \right)e^{-\varkappa \tau} \widehat F_{ 0}(\tau) d\tau.
  \end{equation}
  holds.     
  Using a similar procedure we used to obtain inequalities \eqref{eq:estimate:7},
  and \eqref{eq:estimate:8} respectively we conclude
  \begin{equation}
    |\widehat u_0(t)|\leq \frac{2}{2\varkappa^2}
    \sup_{[0,\infty)}\|a(t,\cdot)\|_{L^\infty}C(2\alpha\beta).
  \end{equation}
  Thus by the smallness condition \eqref{eq:cond:1}, we obtain
  \begin{equation}
    |\widehat u_0(t)|\leq\frac{\alpha\beta\epsilon_1}{3}. 
  \end{equation}
  and observe that inequality \eqref{eq:decay} holds with
  $\widetilde \epsilon=\frac{\alpha\beta\epsilon_1}{3}$.

\end{proof}

\biblio

%% file: section5-blowup.tex
\section{Blowup of solutions even for small initial data}
\label{sec:blow-solut-large}
In the previous sections, \ref{subsec:homogeneous} and
\ref{sec:glob-exist-bootstr}, we proved the global existence (and
uniqueness) of classical solutions in the Sobolev spaces $H^m$.
It is important to emphasize, that the smallness of $a(t,x)$ played an
essential role in the proof. 
That is why we want to drop to the smallness assumption on $a(t,x)$
and investigate its consequence. 
Our main result can be stated as follows.

\begin{thm}[Blowup in finite time]
  \label{thm:blow-up}
  Let $u$ be the solution to the Cauchy problem
  \eqref{eq:wave:1}--\eqref{eq:wave:2} in the interval $[0,T)$, where
  $0<T\leq \infty$ and assume the following conditions:
  \begin{equation}
    \label{eq:section3A-blowup:8}
    0<a_{0}\leq a(t,x),  \qquad \forall   t\geq 0,
  \end{equation}
  \begin{equation}
    \label{eq:blowup:3}
    1+f(x)>0, \quad \Delta f(x) \geq 0, \qquad x\in \setT^3,
  \end{equation}
  \begin{equation}
   \varkappa(1+f(x))+g(x)\geq 0,
  \end{equation} 
  \begin{equation}
    \label{eq:section3A-blowup:7}
    \widehat g_{0}> 0
  \end{equation}
  and
  \begin{equation}
    \label{eq:section3A-blowup:9}
    \widehat g_0^2-\frac{a_0}{2}(1+\widehat f_0)^4\leq 0.
  \end{equation}
  Then for sufficiently large $a_0$, $T$ is finite, and moreover the
  following holds:
  \begin{equation}
    \label{eq:section6-blowup:1}
    \lim_{t\uparrow T}\|u(t,\cdot)\|_{H^{m+1}}=\infty.
  \end{equation}
\end{thm}
We recall that $\widehat f_0$ and $\widehat g_0$ are the zero Fourier
coefficients of $f$ and $g$ respectively. 
Also, note that condition \eqref{eq:blowup:3} implies that
$1+\widehat f_0>0$.
\begin{rem}
\label{rem:section5-blowup:1} \hfill
  \begin{enumerate}[label=\arabic*.]
    \item Note that for large $a(t,x)$ blow-up occurs in finite time
    even if the initial data are small. 
    \item We actually can neglect condition
    \eqref{eq:section3A-blowup:9} since most likely it holds when
    $a_0$ is large.
    \item

    We proved the blow up under the assumptions that $\Delta f(x)\geq 0$ and 
    $\widehat g_0>0$, where
    $\partial_t u(0,x)=g(x) $. 
    Our conjecture is that the blow up holds even without those restrictions.
  \end{enumerate}
\end{rem}

\textsc{Proof Sketch:}
By the local existence theorem, Theorem \ref{thr:section3-fourier:1},
there exists a regular unique solution $u$ to the Cauchy problem,
\eqref{eq:wave:1}--\eqref{eq:wave:2}, namely
$u(t,\cdot)\in L^\infty([0,T]; H^{m+1}(\setT^3)\cap    C^{0,1}([0,T]; 
H^{m}(\setT^3)$

We adopt an idea of Yagdjian \cite{Yagdjian_05} that was used for
a different wave equation, and we set
\begin{equation}
  \label{eq:blow:1}
  F(t)\upperrel{\rm def}{=}\frac{1}{(2\pi)^3}\int_{\setT^3} u(t,x)dx=\widehat 
u_{ 0}(t)
\end{equation}
and we shall derive differential inequality \eqref{eq:blow:6-a} for
$F$.
For that purpose, we need to apply Jensen's inequality to the right
hand side of equation \eqref{eq:blow:3}. 
In order to do so, we have to ensure that the term $1+u(t,x)$ is
nonnegative in the existence interval $[0,T)$ and we prove it in
Section \ref{sec:positivity}.

Finally, we use Lemma \ref{lem:4} below, which states that a function
that satisfies differential inequality \eqref{eq:blow:6-a} blows up in
finite time. Since
\begin{equation}
\label{eq:blowup:9}
 \|u(t)\|_{H^{m+1}}^2=|F(t)|^2+\sum_{0\neq k\in\setZ^3}|k|^{2(m+1)} | \widehat 
u_k(t)|^2\geq |F(t)|^2
\end{equation}
holds, the blow up of $F$ implies the blow up of the solution to the
Cauchy problem \eqref{eq:wave:1}--\eqref{eq:wave:2}.

\begin{proof}[Proof of Theorem \ref{thm:blow-up}]
  We start to derive a differential inequality for $F$ that is defined by
  equation \eqref{eq:blow:1}. 
  First, note that
  \begin{equation}
    F'(t)=\frac{1}{(2\pi)^3}\int_{\setT^3} \partial_t u(t,x)dx =\widehat u_{ 
      0}'(t) \quad \text{ and}\quad 
    F''(t)=\frac{1}{(2\pi)^3}\int_{\setT^3} \partial_{tt} u(t,x)dx.
  \end{equation}
  It is a well known fact that the integral of the Laplacian over the
  $\setT^3$ is zero, or in other words,
  \begin{equation}
    \label{eq:section5-blowup1}
    \int_{\setT^3} \Delta u(t,x)dx =0.
  \end{equation}
  One way to prove equation \eqref{eq:section5-blowup1} is to expand $u$ to its
  Fourier series \eqref{eq:fs:1}, then we observe that the zero
  coefficient $\Delta u $ is zero (another possibility to prove
  equation \eqref{eq:section5-blowup1} is to use Gauss' theorem).
  So we obtain 
  \begin{equation}
    \label{eq:blow:3}
    F''+2\varkappa F'=\frac{1}{(2\pi)^3}\int_{\setT^3}\left(\partial_{tt} 
      u+2\varkappa
      \partial_t 
      u-\Delta u\right)dx
    =\mathord{e^{-\varkappa t}}\frac{1}{(2\pi)^3} 
    \int_{\setT^3}a(t,x)(1+u(t,x))^3dx.
  \end{equation}
  Our idea is to estimate the right hand side of equation
  \eqref{eq:blow:3} by Jensen's inequality (see
  e.~g.~\cite[Ch.~2]{Lieb-Loss}), with the convex function $s^3$,
  where $s=\sqrt[3]{a(t,x)}(1+u(t,x))$. 
  The function $s^3$ is convex, only for $s\geq 0$.
  Recall that $0<a(t,x)$ holds by assumption
  \eqref{eq:section3A-blowup:8}. 
  The proof of the positivity of $1+u(t,x)$ is more involved and we refer to
  Theorem \ref{thm:positivity} in  Section \ref{sec:positivity}.

   Applying Jensens's inequality as outlined above we  obtain
  \begin{equation}
    \label{eq:section5-blowup3}
    \begin{split}
      & \frac{1}{(2\pi)^3}\int_{\setT^3} a(t,x)\left(1+u(t,x)\right)^3 dx
      \geq \left( \frac{1}{(2\pi)^3}\int_{\setT^3}
        a^{1/3}(t,x)\left(1+u(t,x)\right)dx\right)^3\\
      & \geq a_0 \left(\frac{1}{(2\pi)^3}\int_{\setT^3}
        \left(1+u(t,x)\right)dx\right)^3=a_0\left(1+F(t)\right)^3.
    \end{split}
  \end{equation}

  Hence, we have obtained the following initial differential
  inequality
  \begin{subequations}
    \begin{align}
      \label{eq:blow:6-a}
      & F''+2\varkappa F'\geq e^{-\varkappa t} a_0\left(1+F\right)^3
      \\
      \label{eq:blow:6-b}
      & 
        F(0)=\widehat f_0, \quad \ F'(0)=\widehat g_0
    \end{align}
  \end{subequations}
  for $t\in [0,T)$. 
  We now use Lemma \ref{lem:4}, stating that $F$ blows up in a finite
  time interval, then $u(t,x)$ also blows up using equation \eqref{eq:blowup:9}.

\end{proof}

It remains to state and  to prove Lemma \ref{lem:4}.

\begin{lem}[Blow up for the associated differential inequality]
  \label{lem:4}
  Let $F$ satisfy the differential inequality \eqref{eq:blow:6-a} in
  the interval $[0,T)$, where $0<T\leq \infty$ and with the initial data
  \eqref{eq:blow:6-b}. 
  Suppose assumptions \eqref{eq:blowup:3}--\eqref{eq:section3A-blowup:9} of Theorem \ref{thm:blow-up} hold,
  then for sufficiently large $a_0$, $T$ is finite, and moreover
  \begin{equation}
    \label{eq:blow:6}
    \lim_{t \uparrow T}F(t)=\infty.
  \end{equation}
\end{lem}

Note, that the differential inequality \eqref{eq:blow:6-a} contains a
strong damping term $e^{-\varkappa t}$. 
So we want to know whether this damping term prevents the blowup in
finite time. 
It turns out, however, that this is not the case.

\textsc{Proof sketch:} The proof consists of three main steps. 
In the first step, we show that $F'(t)$ is a positive function in the
interval of existence $[0,T)$.
In the second step, we shall make a variable change in order to
transform the differential inequality \eqref{eq:blow:6-a} to an
inequality without a first-order term. 
This together with the positivity of $F'$ will enable us in the third
step to integrate the inequalities and to estimate the time of the
blow up.
\begin{proof}[Proof of Lemma \ref{lem:4}]
\hfill
\begin{enumerate}[label=\textbf{Step \arabic*},wide, itemsep=4ex, labelwidth=!, labelindent=0pt]
  \item \label{item:section5-blowup:1}
  We claim that $F^\prime(t)>0$ holds in
  the existence interval $[0,T)$. 
  To see that we set
  \begin{equation}
    \label{eq:section5-blowup:1}
    T^\ast=\sup\{T_1: F'(t)>0 \quad \text{for}\ \ t\in [0,T_1), 0\leq T_1\leq T\}.
  \end{equation}
  By the assumptions of Lemma \ref{lem:4},
  $F^\prime(0)=\widehat g_0>0$, hence $T^\ast>0$. 
  We now assume by contradiction that $T^\ast<T$, then by the
  continuity of $F^\prime(t)$, we conclude  that $F^\prime(T^\ast)=0$. 
  Recall that by assumption \eqref{eq:blowup:3},
  $1+F(0)=1+\widehat f_0>0$, hence $1+F(T^\ast)\geq 1+F(0)>0$ and
  consequently
  \begin{equation}
    \label{eq:section5-blowup:3} F^{\prime\prime}(T^{\ast}) =
    F^{\prime\prime}(T^{\ast}) +2\varkappa F'(T^{\ast})\geq e^{-\varkappa 
      T^{\ast}}
    a_0\left(1+F(T^{\ast})\right)^3>0.
  \end{equation}
  This implies that $F$ attains a local minimum at time $T^\ast$.
  But this is impossible since $F$ is an increasing function. 
  Therefore we conclude that $T^\ast =T$.

  \item\label{item:section5-blowup:2} We start with a variable change
  of the form $\tau=\omega(t)$ and define a function $G(\tau)$ such that
  $F(t)=G(\omega(t))$. 
  Then $F^\prime (t)= \frac{ dG}{d \tau}(\omega(t))\omega'(t)$ and
  $F^{\prime\prime} (t)= \frac{ d^2G}{d \tau^2}(\omega(t))(\omega'(t))^2+
  \frac{ dG}{d \tau}(\omega(t))\omega^{\prime\prime}(t)$, which leads to
  \begin{equation}
    \label{eq:blow:2}
    F^{\prime\prime}+2\varkappa F^\prime= \frac{d^2 G}{d\tau^2}\left(\omega^\prime 
    \right)^2 
    +\frac{d G}{d\tau}\left(\omega^{\prime\prime}+2\varkappa\omega^\prime\right).
  \end{equation}
  Now we choose $\omega $ such that it satisfies the equation
  $\omega^{\prime\prime}+2\varkappa\omega^\prime=0$, and in order to
  obtain an one-to-one transformation for $t >0$ we require that
  $\omega'>0$.
  It is straightforward to calculate its general solution
  \begin{equation*}
    \omega(t)=C_1e^{-2\varkappa t}+C_0,
  \end{equation*}
  for which we choose $C_1=-1 $ and $ C_0=2$, and then
  \begin{equation}
    \label{eq:section3A-blowup:5}
    \tau=\omega(t)=-e^{-2\varkappa t}+2.
  \end{equation}
  Note that $\omega $ maps $[0,\infty) $ onto $[1,2)$ in a one-to-one
  manner. 
  Taking into account equation \eqref{eq:blow:2}, inequality
  \eqref{eq:blow:6-a} is equivalent to
  \begin{equation}
    \label{eq:blow:8a}
    \frac{d^2  G}{d\tau^2}\left(\omega^\prime\right)^2\geq e^{-\varkappa t} 
    a_0(1+G)^3
  \end{equation}
  or
  \begin{equation}
    \label{eq:blow:9}
    \frac{d^2  G}{d\tau^2}\geq \frac{e^{-\varkappa t}}{(2\varkappa e^{-2\varkappa 
        t})^2} 
    a_0(1+G)^3=\frac{e^{3\varkappa
        t}}{4\varkappa^2}a_0(1+G)^3
    \geq\frac{     a_0}{4\varkappa^2}(1+G)^3.
  \end{equation}
  In order to simplify the notation, we set
  $\frac{d G}{d\tau}=G^\prime$ and
  $G^{\prime\prime}=\frac{d^2 G}{d\tau^2}$, then $G$ satisfies the
  initial value inequality
  \begin{subequations}
    \begin{align}
      \label{eq:blow:5}
      & G''\geq  \frac{a_0}{4\varkappa^2}(1+G)^3 \\
      \label{eq:blow:4}
      & G(1)=\widehat f_0, \quad G'(1)=\frac{\widehat g_0}{2\varkappa}.
    \end{align}
  \end{subequations}
  We are now in a position to show the blow up for $G$ at some
  $1<\tau_0<2$, and consequently, $F$ will blow up at
  $t_0=\omega^{-1}(\tau_0)$.
  
   \vskip 5mm
  \item\label{item:section5-blowup:3} We will show that if $a_0$ is
  sufficiently large, then there exits $\tau_0<2$ such that
  \begin{equation}
    \label{eq:blow:7}
    \lim_{\tau \uparrow \tau_0} G(\tau)=\infty.
  \end{equation}
  Note that $F'>0$ holds in the existence interval as it was proven in
  \ref{item:section5-blowup:1}, and since $F'=G'\frac{d\omega}{d t}$ and
  $\omega^\prime>0$, we can conclude that $G'>0$ holds. 
  Thus we can multiply inequality \eqref{eq:blow:5} by $G'$ and obtain
  \begin{equation*}
    G'' G'\geq \frac{a_0}{4\varkappa^2}(1+G)^3 G'.
  \end{equation*}
  Integrating both sides of inequality \eqref{eq:blow:5} from $1 $ to $\tau$, we
  conclude that
  \begin{equation*}
    \frac{1}{2}\left(\left(G'(\tau)\right)^2-\left(G'(1)\right)^2\right)\geq
    \frac{a_0}{16 
      \varkappa^2}\left(\left(1+G(\tau)\right)^4-\left(1+G(1)\right)^4\right).
  \end{equation*}
  Taking into account the initial values \eqref{eq:blow:4} we observe
  that
  \begin{equation}
    \begin{split}
      \label{eq:blow:14}
      \left(G'(\tau)\right)^2 &\geq
      \frac{a_0}{8\varkappa^2}\left(1+G(\tau)\right)^4+\left(\frac{\widehat
          g_0}{2\varkappa
        }\right)^2 -\frac{a_0}{8\varkappa^2}\left(1+\widehat f_0\right)^4\\
      & =\frac{a_0}{8\varkappa^2} \left\{\left(1+G(\tau)\right)^4 -
        \left\{ \left(1+\widehat f_0\right)^4-\frac{2\widehat
            g_0^2}{a_0}\right\}\right\}.
    \end{split}
  \end{equation}
    
  The expression
  $(1+\widehat f_0)^4-\frac{2\widehat g_0^2}{a_0}\geq 0$ by assumption
  \eqref{eq:section3A-blowup:9}, and in order to simplify the
  calculations we set
  \begin{equation*}
    \lambda^4\eqdef(1+\widehat f_0)^4-\frac{2\widehat g_0^2}{a_0}.
  \end{equation*}
  Now, since $G'(\tau)>0$,
  \begin{equation}
    (1+G(\tau))^4\geq (1+G(1))^4=(1+\widehat f_0)^4\geq \lambda^4,
  \end{equation}
  hence,
  \begin{equation}
    G'(\tau)\geq 
    \frac{\sqrt{a_0}}{\sqrt{8}\varkappa}\left\{
      \left(1+G(\tau)\right)^4-\lambda^4\right\}^{\frac{1}{2}}\geq
    \frac{\sqrt{a_0}}{\sqrt{8}\varkappa}\left\{
      \left(1+G(\tau)\right)^2-\lambda^2\right\},
  \end{equation}
  or
  \begin{equation}  
    \label{eq:blow:12}
    \frac{\sqrt{8}\varkappa}{\sqrt{a_0}}\dfrac{G'(\tau)}{
      \left(1+G(\tau)\right)^2-\lambda^2}=\frac{\sqrt{8}\varkappa}{\sqrt{a_0}}\frac{
      G'(\tau)}{2\lambda}\left(\dfrac{1}{
        \left(1+G(\tau)\right)-\lambda}-\dfrac{1}{
        \left(1+G(\tau)\right)+\lambda}\right)
    \geq 1.
  \end{equation}
  Integration of both sides of equation \eqref{eq:blow:12} results in
  \begin{equation}
    \frac{\sqrt{2}\varkappa}{\lambda\sqrt{a_0}}\left\{\ln\left(\frac{
          1+G(\tau)-\lambda}{ 1+G(\tau)+\lambda}\right)-\ln\left(\frac{
          1+\widehat f_0-\lambda}{ 1+\widehat f_0+\lambda}\right)\right\}\geq \tau-1,
  \end{equation}
  or
  \begin{equation}
    \label{eq:blow:11}
    \ln\left(\frac{
        1+G(\tau)-\lambda}{ 1+G(\tau)+\lambda}\right)\geq 
    \frac{\lambda\sqrt{a_0}}{\sqrt{2}\varkappa}(\tau-1)+\ln\left(\frac{
        1+\widehat f_0-\lambda}{ 1+\widehat f_0+\lambda}\right).
  \end{equation}
  Note that by inequality \eqref{eq:section3A-blowup:7}
  $1+\widehat f_0-\lambda>0$ so inequality \eqref{eq:blow:11} is well defined.
  We set
  \begin{equation}
    \beta\eqdef\frac{1+\widehat f_0-{\lambda}}{1+\widehat 
      f_0+{\lambda} },
  \end{equation}
  then inequality \eqref{eq:blow:11} is equivalent to
  \begin{equation}
    \frac{1+G(\tau)-\lambda}{ 1+G(\tau)+\lambda}\geq 
    e^{\frac{\lambda\sqrt{a_0}}{\sqrt{2}\varkappa}(\tau-1)}\beta,
  \end{equation}
  and from this inequality we get that
  \begin{equation}
    \label{eq:blow:10}
    {1+G(\tau)}\geq \dfrac{\lambda \left( 
        e^{\frac{\lambda\sqrt{a_0}}{\sqrt{2}\varkappa}(\tau-1)}\beta+1\right)}{1-\beta 
      e^{\frac{\lambda\sqrt{a_0}}{\sqrt{2}\varkappa}(\tau-1)}}.
  \end{equation}
  The right hand side of inequality \eqref{eq:blow:10} blows up at
  time $\tau_0$ for which $\ln(1/\beta)= \frac{\lambda\sqrt{a_0}}{\sqrt{2}\varkappa}(\tau_0-1)$ or
  \begin{equation}
    \tau_0=\frac{\sqrt{2}\varkappa}{\lambda\sqrt{a_0}}\ln\left(\frac{1}{\beta}
    \right)+1.
  \end{equation}
  Thus blow up will occur at time $\tau_0$, however, in order that it
  will be finite in $t$, we have to assure that $ \tau_0<2$, that is,
  \begin{equation}
    \label{eq:blow:8}
    \frac{\sqrt{2}\varkappa}{\lambda\sqrt{a_0}}\ln\left(\frac{1}{\beta}
    \right)<1.
  \end{equation}

  So we now estimate this expression. 
  Recall that
  \begin{equation}
    \lambda^4= (1+\widehat f_0)^4-\frac{2\widehat g_0^2}{a_0}=(1+\widehat 
    f_0)^4\left(1-\frac{2\widehat g_0^2}{(1+\widehat f_0)^4a_0}\right).
  \end{equation}
  We set
  \begin{equation}
    \label{eq:blow:13}
    z=\frac{2\widehat g_0^2}{(1+\widehat f_0)^4a_0}
  \end{equation}
  and we note that $z$ becomes smaller as $a_0$ grows. 
  Hence
  \begin{equation}
    \lambda =(1+\widehat f_0)(1-z)^{\frac{1}{4}}=(1+\widehat 
    f_0)\left(1-\frac{1}{4}z+o(z)\right),
  \end{equation}
  and
  \begin{equation}
    \frac{1}{\beta}  =\frac{(1+\widehat f_0)+{\lambda}}{(1+\widehat 
      f_0)-{\lambda} }=\frac{1+(1-z)^{\frac{1}{4}}}{1-(1-z)^{\frac{1}{4}}}
    =\frac{2-\frac{z}{4}+o(z)}{\frac{z}{4}+o(z)}=\frac{8}{z}-1+o(z).
  \end{equation}
  We now express $a_0$ by $z$ through condition \eqref{eq:blow:13},
  then
  \begin{equation}
    \label{eq:blow:21}
    \frac{\sqrt{2}\varkappa}{\lambda\sqrt{a_0}}\ln\left(\frac{1}{\beta}
    \right)=\frac{\varkappa (1+\widehat f_0)\sqrt{z}}{(1-z)^{\frac{1}{4}}\widehat 
      g_0}\left( 
      \ln\left(\frac{8}{z}-1+o(z)\right)\right) 
    \to 0, \quad \text{as} \ z\to 0.
  \end{equation}
  Consequently, inequality \eqref{eq:blow:8} holds for $a_0$ sufficiently large
  and there exists $\tau_0<2$ such that equation \eqref{eq:blow:7} is true.
\end{enumerate}
\end{proof}

\biblio

%% file: section6-positivity.tex
\section{Positivity  of the scalar field function}
\label{sec:positivity}
The proof of the blowup of classical solution presented in Section
\ref{sec:blow-solut-large} depends on Jensen inequality. 
This inequality states that
\begin{equation}
  \label{eq:section6-positivity:1}
 \Phi\left(\int f  d \mu\right)\leq \int \Phi\left(f \right) d \mu,
\end{equation}
holds, whenever $\Phi$ is a convex function and $\mu$ is a probabilistic
measure. 
In our case, the function $\Phi(s)=s^3$ is a convex function only if
$ s\geq 0$. 
Hence, in order to apply \eqref{eq:section6-positivity:1}, we need to
show that $ 1+u(t,x)\geq 0$.
Since the proof of this inequality is a bit lengthy and requires
additional tools, we have moved the proof to a separate section.  

Another important issue with the positivity of $1+u(t,x)\geq 0$ is the
following.
The metric of the spacetime is given by
$g_{\alpha\beta}=\phi^2\eta_{\alpha\beta}$, where $\phi(t,x)=e^{\varkappa t}(1+u(t,x))$. 
Hence, if $1+u(t,x)=0$, then metric vanishes and the solution in that
case has no physical meaning.

We recall that $u $ satisfies the initial value problem
\eqref{eq:wave:1}--\eqref{eq:wave:2}, and the known existence theorem
for semi-linear wave equations, \cite[Theorem 6.4.11]{Hormander_1997},
assures the existence and uniqueness of $C^2 $-solution in a certain
time interval $[0, T)$. 
Our aim is to show that $ 1+u(t,x)\geq 0$ provided that $a(t,x)>0 $ and
the initial data $ 1+u(0,x)=1+f(x)$ are positive. 
However, the condition $1+f(x)>0$ is not sufficient and further
conditions are needed. 
The question of which additional conditions to impose has also been
discussed in other publications in which the  solution of the
linearized equation is required to be positive, see for example
\cite{Caffarelli_Friedman_86} and the celebrated paper by John
\cite{John_79}. 

In the following we will present the  main result of this section.
Since our major interest  here is the positivity, we assume that initial  
data $ f $ and $ g $  are sufficiently smooth in this section. 

\begin{thm}[The positivity of $1+u(t,x)$]
  \label{thm:positivity}
  Assume $ u(t,x)$ is a unique $C^2 $ solution of initial value problem
  \eqref{eq:wave:1}--\eqref{eq:wave:2} for $ t\in [0,T)$,
  $x\in \setT^3$ and for some positive $T$. 
  Moreover, assume that 
  \begin{align}
    \label{eq:section6-positivity:6}
    a(t,x)&>0 \\
    \label{eq:positivity:5}
    1+f(x)&>0, \quad g(x) \geq 0, \quad x\in \setT^3\\
  \label{eq:positivity:6}
   \Delta f(x)&\geq 0, \quad x\in \setT^3.
  \end{align}
  Then
  \begin{equation}
    \label{eq:section6-positivity:4}
    1+u(t,x)>0 \quad \text{ for } \ (t,x)\in [0,T)\times \setT^3.
  \end{equation}
\end{thm}
The proof of Theorem \ref{thm:positivity} is based on the and the
application of Kirchhoff's formula for linear wave equations in
$\setR^3$ \cite[\S 5]{john86:_partial} and this is why we have to
linearize equation \eqref{eq:wave:1}. 
We first transform the equation \eqref{eq:wave:1} with the damping
term, to an appropriate form by setting $\phi(t,x)=e^{\varkappa t}(1+ u(t,x))$. 
We recall that $\phi $ satisfies the field equation
\eqref{eq:Nordstrom:3} with the cosmological constant $\varkappa^2$, thus the
initial value problem \eqref{eq:wave:1}--\eqref{eq:wave:2} is
equivalent to
\begin{subequations}
  \begin{align}
    \label{eq:positivity:1}
    \phi_{t t}-\Delta \phi&=e^{-3 \varkappa t}  a(t, x) \phi^{3}+\varkappa^{2} \phi \\
    \label{eq:positivity:18}
    \phi(0, x)&=1+f(x)    \\
    \label{eq:section6-positivity:2}
    \phi_{t}(0, x)&=h(x),
  \end{align}
\end{subequations}
where
\begin{equation}
 h(x)\upperrel{\text{\rm def}}{=}\varkappa(1+f(x))+g(x).
\end{equation} 

The linearization of \eqref{eq:positivity:1}, \eqref{eq:positivity:18} and \eqref{eq:section6-positivity:2} results in
the system 
\begin{subequations}
\begin{align}
    \label{eq:positivity:2}
    v_{tt} - \Delta v&=P(t, x) \\
    \label{eq:positivity:22}
    v(0, t)&=1+f(x)\\
  \label{eq:section6-positivity:3}
   v_{t}(0, x)&=h(x),
  \end{align}
\end{subequations}
where $P(t,x)$ denotes the linearization of the nonlinear right hand
side of \eqref{eq:positivity:1} and whose precise form it of no
importance. 
By Kirchhoff's formula (see e.g.~\cite[\S 5]{John_79}) solution of the
above system is given by
\begin{equation}
  \label{eq:positivity:4}
  \begin{aligned}
    v(t, x)=& \frac{t}{4 \pi} \int_{|\xi|=1} h(x+t \xi) d \omega_\xi+\frac{\partial}{\partial t}\left(\frac{t}{4 \pi} \int_{|\xi|=1}(1+f(x+t\xi)) d \omega_\xi\right)\\
    &+\frac{1}{4 \pi} \int_{0}^{t}(t-s)\left(\int_{|\xi|=1} P(s, x+(t-s) \xi) d
      \omega_\xi \right) d s,
  \end{aligned}
\end{equation}
where $ d \omega_\xi$ is the Lebesgue  measure of the unite sphere $\mathbb{S}^2$.

\begin{rem}
  In this paper, however, we deal with spatially periodic solutions,
  while Kirchhoff's formula \eqref{eq:positivity:4} provides solutions
  in $\setR^3$. 
  But if initial data \eqref{eq:positivity:22} and
  \eqref{eq:section6-positivity:3}, and the right hand side of
  \eqref{eq:positivity:2} are periodic, then it follows from
  \eqref{eq:positivity:4} that $v(t,x) $ is also a periodic function
  of the space variable $x$. 
\end{rem}
\textsc{Proof sketch:} Our proof's strategy can be described as follows:
\begin{enumerate}[label=\arabic*.]
  \item We first show that the solution to the homogeneous initial
  value problem
  \eqref{eq:positivity:2} --~\eqref{eq:section6-positivity:3} is
  positive, if the conditions \eqref{eq:positivity:5} and 
  \eqref{eq:positivity:6} are met.
    
  \item We then construct a monotone sequence,
  $0<\phi_n\leq \phi_{n+1}$, by an iteration of the linearized equation. 
  
  \item We show that whenever $\phi(x,t)>0$, then
  \begin{equation}
    \label{eq:skech:1}
    \phi_n(t,x)\leq \phi(t,x).
  \end{equation} 
  We then  use  \eqref{eq:skech:1}  to show that
  $\phi(t,x)>0 $ in the entire existence interval.  
  \item \label{item:section6-positivity:3} Finally, the positivity of
  $\phi$ implies that $1+u(t,x)=e^{-\varkappa t}\phi(t,x)>0$.
\end{enumerate}
  
\subsection*{The iteration scheme}
We denote the right hand side of \eqref{eq:positivity:1} by $G(\phi)$, that is, 
\begin{equation*}
  G(\phi)=G(\phi, t, x)=e^{-3 k t} a(t, x) \phi^{3}+k^{2} \phi.
\end{equation*}
Note that $G(0)=0$ and 
  \begin{equation*}
\frac{\partial}{\partial \phi} G(\phi)=3 e^{-3 k t} a(t, x) \phi^{2}+k^{2}.
\end{equation*}
Hence,  if $a(t,x)>0 $, then $ G $ is an  increasing function of $\phi$ and  non-negative for $\phi\geq 0 $.   
We  now let $\phi_0(t,x)=0$ and set $\phi_{n+1}$ to be the solution to the linear equation
\begin{equation}
  \label{eq:positivity:8}
  \begin{cases}
     & \left(\phi_{n+1}\right)_{t t}-\Delta \phi_{n+1} =G\left(\phi_{n}\right) \\
     & \phi_{n+1}(0, x)                      =1+f(x)          \\
     & (\phi_{n+1})_{t}(0, x)                 =h(x)
    \end{cases}.
  \end{equation}
  
\subsubsection*{\textsc{Step 1: The free wave equation.}}
\label{sec:step-1:-free}

The first step of the iteration  consists of showing that the free wave equation 
\begin{equation}
  \label{eq:positivity:30}
  \begin{cases}
     & \left(\phi_{1}\right)_{t t}-\Delta \phi_{1} =G\left(0\right)=0 \\
     & \phi_{1}(0, x)                      =1+f(x)            \\
     & (\phi_{1})_{t}(0, x)                 =h(x)
    \end{cases}
  \end{equation}
has a positive solution.
\begin{prop}
  \label{prop:6.1}
Under conditions \eqref{eq:positivity:5} and \eqref{eq:positivity:6} it follows that $\phi_1(t,x)>0$.
\end{prop}
\begin{proof}
  Taking the time derivative in  Krichhoff's formula \eqref{eq:positivity:4}, we observe  that
  \begin{equation}
\label{eq:section6-positivity:5}
    \phi_1(t, x) =\frac{t}{4 \pi} \int_{|\xi|=1} \left(h(x+t \xi)+ \xi\cdot \nabla f(x+t\xi)\right) d \omega_\xi+\frac{1}{4 \pi} 
    \int_{|\xi|=1}(1+f(x+ t\xi) )d \omega_\xi.
  \end{equation}
  Applying   now  Gauss Divergence Theorem to  term $ \xi\cdot \nabla f(x+t\xi)$, we obtain that
   \begin{equation}
\label{}
    \phi_1(t, x) =\frac{t}{4 \pi} \int_{|\xi|=1} \left(h(x+t \xi)\right) d \omega_\xi
    + \frac{t^2}{4 \pi} \int_{|\xi|\leq1} \Delta f(x+t \xi) d \xi
    +\frac{1}{4 \pi} \int_{|\xi|=1}(1+f(x+ t\xi) )d \omega_\xi.
  \end{equation}
Hence   conditions \eqref{eq:positivity:5} and \eqref{eq:positivity:6} imply that 
$\phi_1(t,x)>0$.
\end{proof}

\subsubsection*{\textsc{Step 2: Monotonicity.}}
\label{sec:textscst-2:-monot}
\begin{prop}
 \label{prop:monoton}
 Assume  conditions \eqref{eq:section6-positivity:6}, \eqref{eq:positivity:5} and \eqref{eq:positivity:6} are satisfied 
 then the sequence $\{\phi_n\}$   defined by \eqref{eq:positivity:8} is
 monotone, in other words
 \begin{equation*}
  \phi_n(t,x)\leq \phi_{n+1}(t,x).
 \end{equation*}
holds
\end{prop}
\begin{proof}
  The proof is obviously done by induction. 
  We already have proved that $\phi_1(t,x)>\phi_0(t,x)\equiv 0$ holds.
  Assume $\phi_{n-1}\leq \phi_n$, then $G(\phi_{n-1})\leq G(\phi_n)$ since
  both functions are positive and  $G $ is increasing. 
  Hence, by Kirchhoff's formula \eqref{eq:positivity:4}, we observe  that
  \begin{equation}
    \begin{aligned}
      \phi_{n+1}(t, x)=& \frac{t}{4 \pi} \int_{|\xi|=1} h(x+t \xi) d \omega_\xi+\frac{\partial}{\partial t}\left(\frac{t}{4 \pi} \int_{|\xi|=1}(1+f(x+t\xi)) d \omega_\xi\right)\\
      &+\frac{1}{4 \pi} \int_{0}^{t}(t-s)\left(\int_{|\xi|=1} G(\phi_{n})(s,
        x+(t-s) \xi) d \omega_\xi \right) d s \\ & \geq
      \frac{t}{4 \pi} \int_{|\xi|=1} h(x+t \xi) d \omega_\xi+\frac{\partial}{\partial t}\left(\frac{t}{4 \pi} \int_{|\xi|=1}(1+f(x+t\xi)) d \omega_\xi\right)\\
      &+\frac{1}{4 \pi} \int_{0}^{t}(t-s)\left(\int_{|\xi|=1} G(\phi_{n-1})(s,
        x+(t-s) \xi) d \omega_\xi \right) d s \\ & =\phi_n(t,x).
    \end{aligned}
  \end{equation}
which finishes the proof.
\end{proof}
\textsc{Step 3:} It remains to show that  $\phi(t,x)>0$ holds.
\begin{prop}
  \label{porp:positive:2}
  Assume $f(x)$,  $h(x) $ and $a(t,x)$  are periodic functions that satisfies the  assumptions of 
  Proposition \ref{prop:monoton}. 
  Let $\phi $ be a $C^2$ solution to initial value problem 
  \eqref{eq:positivity:1}-\eqref{eq:section6-positivity:2} in the interval
  $[0,T)$. 
  Then
  \begin{equation}
    \phi(t,x)>0, \quad t\in [0,T) \ \text {and all}\ x\in\setT^3. 
  \end{equation}
\end{prop}
\begin{proof}
  We set
  \begin{equation}
    \label{eq:positivity:20}
    T^{*}=\sup \left\{T_{1}: \phi(t, x)>0\right.
    \mbox{for} \left.t \in\left[0, T_{1}\right) \ \text{and } \forall x \in \setT^{3}\right\}.
  \end{equation}
    
  Since $1+f(x)>0$, then by continuity, $0<T^\ast \leq T$. 
  The proof consists essentially of two steps. 
  In the first one we show that $\phi_n(t,x)\leq \phi(t,x)$ for
  $t\in [0,T^\ast)$, and in the second one we show that $T^\ast=T$. 
    
  We prove the first step, again, by induction. 
  Proposition \ref{prop:6.1} implies that  $0=\phi_0(t,x)<\phi(t,x) $ for $ t\in [0,T^\ast)$. 
  Now, assume $\phi_{n-1}(t,x)\leq \phi(t,x)$ for $ t\in [0,T^\ast)$, then in a
  similar way to the proof of the monotonicity, Proposition
  \ref{prop:monoton}, we obtain that
  \begin{equation}
    \begin{aligned}
      \phi_{n}(t, x)=& \frac{t}{4 \pi} \int_{|\xi|=1} h(x+t \xi) d \omega_\xi+\frac{\partial}{\partial t}\left(\frac{t}{4 \pi} \int_{|\xi|=1}(1+f(x+t\xi)) d \omega_\xi\right)\\
      &+\frac{1}{4 \pi} \int_{0}^{t}(t-s)\left(\int_{|\xi|=1} G(\phi_{n-1})(s,
        x+(t-s) \xi) d \omega_\xi \right) d s \\ & \leq
      \frac{t}{4 \pi} \int_{|\xi|=1} h(x+t \xi) d \omega_\xi+\frac{\partial}{\partial t}\left(\frac{t}{4 \pi} \int_{|\xi|=1}(1+f(x+t\xi)) d \omega_\xi\right)\\
      &+\frac{1}{4 \pi} \int_{0}^{t}(t-s)\left(\int_{|\xi|=1} G(\phi)(s, x+(t-s)
        \xi) d \omega_\xi \right) d s \\ & =\phi(t,x).
    \end{aligned}
  \end{equation}
  The last equality follows from the fact that $\phi$ satisfies the initial value problem 
  \eqref{eq:positivity:1}--\eqref{eq:section6-positivity:2}.
    
  We turn now to the second step. 
  We claim that $T^\ast=T$. 
  If not, then $T^\ast<T $ and we will derive a contradiction. 
  First, we note that there exists a $x_0\in\setT^3$ such that
  $\phi(T^\ast,x_0)\leq 0$. 
  Since otherwise, by continuity, there exists $T^{\ast\ast}>T^\ast$
  such that $\phi(t,x)>0$ for $t\in [0,T^{\ast\ast})$ and all
  $x\in\setT^3$, and that obviously contradicts the definition of
  $T^\ast$ in \eqref{eq:positivity:20}. 
  Now, using monotonicity and the first step, we conclude that
  \begin{equation}
    \label{eq:positivity:13}
    0 <\lim _{t \rightarrow T^{*-}} \phi_{n}\left(t, x_{0}\right) \leq \lim _{t \rightarrow T^{*-}} \phi\left(t, x_{0}\right)=\phi\left(T^{*}, x_{0}\right) \leq 0
  \end{equation}
  holds, which is the desired contradiction.    
\end{proof}

Step \ref{item:section6-positivity:3} is obvious, and that completes
the proof of Theorem \ref{thm:positivity}.

\biblio

%% file: manuscript.bbl
\providecommand{\bysame}{\leavevmode\hbox to3em{\hrulefill}\thinspace}
\providecommand{\MR}{\relax\ifhmode\unskip\space\fi MR }
\providecommand{\MRhref}[2]{%
  \href{http://www.ams.org/mathscinet-getitem?mr=#1}{#2}
}
\providecommand{\href}[2]{#2}
\begin{thebibliography}{AFCP14}

\bibitem[AFCP14]{Felix_Antonio_Calogero_2014}
J.~A. Alcántara~Felix, S.~Calogero, and S.~Pankavich, \emph{Spatially
  homogeneous solutions of the {V}lasov-{N}ordström-{F}okker-{P}lanck system},
  J. Differential Equations \textbf{257} (2014), no.~10, 3700--3729, URL:
  \url{https://doi.org/10.1016/j.jde.2014.07.006}, \href
  {http://dx.doi.org/10.1016/j.jde.2014.07.006}
  {\path{doi:10.1016/j.jde.2014.07.006}}. \MR{3260238}

\bibitem[BCD11]{Bahouri_2011}
H.~Bahouri, J.~Chemin, and R.~Danchin, \emph{Fourier {A}nalysis and {N}onlinear
  {P}artial {D}ifferential {E}quations}, Springer, Heidelberg London New York,
  2011.

\bibitem[BK14]{BK8}
U.~Brauer and L.~Karp, \emph{Local existence of solutions of self gravitating
  relativistic perfect fluids}, Comm. Math. Phys. \textbf{325} (2014), no.~1,
  105--141, URL: \url{https://doi.org/10.1007/s00220-013-1854-3}. \MR{3182488}

\bibitem[BRR94]{Brauer_Rendall_Reula_94}
U.~Brauer, A.~Rendall, and O.~Reula, \emph{The cosmic no--hair theorem and the
  nonlinear stability of homogeneous {Newtonian} cosmological models},
  Classical and Quantum Gravity \textbf{11} (1994), no.~9, 2283--2296, URL:
  \url{http://lanl.arxiv.org//abs/gr-qc/9403050}.

\bibitem[Cal03]{calogero03:_spher}
S.~Calogero, \emph{Spherically symmetric steady states of galactic dynamics in
  scalar gravity}, Class. Quant. Grav. \textbf{20} (2003), 1729--1742, URL:
  \url{https://doi.org/10.1088/0264-9381/20/9/310}, \href
  {http://arxiv.org/abs/math-ph/0301031} {\path{arXiv:math-ph/0301031}}.

\bibitem[CB09]{Choquet-Bruhat_09}
Y.~Choquet-Bruhat, \emph{General {R}elativity and the {E}instein {E}quations},
  Oxford Mathematical Monographs, Oxford University Press, Oxford, 2009.
  \MR{2473363}

\bibitem[CF86]{Caffarelli_Friedman_86}
L.~A. Caffarelli and A.~Friedman, \emph{The blow-up boundary for nonlinear wave
  equations}, Trans. Amer. Math. Soc. \textbf{297} (1986), no.~1, 223--241,
  URL: \url{https://doi.org/10.2307/2000465}, \href
  {http://dx.doi.org/10.2307/2000465} {\path{doi:10.2307/2000465}}. \MR{849476}

\bibitem[CR03]{Calogero_Rein_2003}
S.~Calogero and G.~Rein, \emph{On classical solutions of the
  {N}ordström-{V}lasov system}, Comm. Partial Differential Equations
  \textbf{28} (2003), no.~11-12, 1863--1885, URL:
  \url{https://doi.org/10.1081/PDE-120025488}, \href
  {http://dx.doi.org/10.1081/PDE-120025488} {\path{doi:10.1081/PDE-120025488}}.
  \MR{2015405}

\bibitem[CR04]{Calogero_Rein_2004}
\bysame, \emph{Global weak solutions to the {N}ordström-{V}lasov system}, J.
  Differential Equations \textbf{204} (2004), no.~2, 323--338, URL:
  \url{https://doi.org/10.1016/j.jde.2004.02.011}, \href
  {http://dx.doi.org/10.1016/j.jde.2004.02.011}
  {\path{doi:10.1016/j.jde.2004.02.011}}. \MR{2085540}

\bibitem[DP85]{Dym_P.Mckean_85}
H.~Dym and H.~P.McKean, \emph{Fourier {S}eries and {I}ntegrals}, Academic
  Press, 1985.

\bibitem[EF14]{einstein14:_nords_gravit_stand_differ}
A.~Einstein and A.~Fokker, \emph{{Die Nordströmische Gravitationstheorie vom
  Standpunkt des absoluten Differentialkalküls}}, Annalen der Physik
  \textbf{44} (1914), 321--328.

\bibitem[ER18]{Ebert_Reissig_18}
M.~Ebert and M.~Reissig, \emph{Methods for {P}artial {D}ifferential
  {E}quations}, Birkhäuser, 2018.

\bibitem[FJS21]{Fajman_Jeremie_Jacques-2021}
D.~Fajman, J.~Joudioux, and J.~Smulevici, \emph{The stability of the
  {M}inkowski space for the {E}instein-{V}lasov system}, Anal. PDE \textbf{14}
  (2021), no.~2, 425--531, URL: \url{https://doi.org/10.2140/apde.2021.14.425},
  \href {http://dx.doi.org/10.2140/apde.2021.14.425}
  {\path{doi:10.2140/apde.2021.14.425}}. \MR{4241806}

\bibitem[Fri86]{Friedrich:1986}
H.~Friedrich, \emph{On the existence of n-geodesically complete or future
  complete solutions of {E}instein's field equations with smooth asymptotic
  structure}, Comm. Math. Phys. \textbf{107} (1986), 587--609.

\bibitem[Hö97]{Hormander_1997}
L.~Hörmander, \emph{Lectures on nonlinear hyperbolic differential equations},
  Mathématiques \& Applications (Berlin) [Mathematics \& Applications],
  vol.~26, Springer-Verlag, Berlin, 1997. \MR{1466700}

\bibitem[Joh79]{John_79}
F.~John, \emph{Blow-up of solutions of nonlinear wave equations in three space
  dimensions}, Manuscripta Math. \textbf{28} (1979), no.~1-3, 235--268, URL:
  \url{https://doi.org/10.1007/BF01647974}, \href
  {http://dx.doi.org/10.1007/BF01647974} {\path{doi:10.1007/BF01647974}}.
  \MR{535704}

\bibitem[Joh86]{john86:_partial}
F.~John, \emph{Partial differential equations}, Springer, 1986.

\bibitem[Kat75]{KATO}
T.~Kato, \emph{{The Cauchy Problem for Quasi--Linear Symmetric Hyperbolic
  Systems}}, Archive for Rational Mechanics and Analysis \textbf{58} (1975),
  181--205.

\bibitem[LL01]{Lieb-Loss}
E.~Lieb and M.~Loss, \emph{Analysis}, 2 ed., Graduated Studies in Mathematics,
  vol.~14, AMS, Providence, Rhode Island, 2001.

\bibitem[Nor13]{nordstrom13:_zur_theor_gravit_stand_relat}
G.~Nordström, \emph{Zur {T}heory der {G}ravitation vom {S}tandpunkt des
  {R}elativitätsprinzip}, Annalen der Physik \textbf{43} (1913), 533--554.

\bibitem[Rau12]{rauch12:_hyper}
J.~Rauch, \emph{Hyperbolic partial differential equations and geometric
  optics}, Graduate Studies in Mathematics, vol. 133, American Mathematical
  Society, Providence, RI, 2012, URL: \url{https://doi.org/10.1090/gsm/133},
  \href {http://dx.doi.org/10.1090/gsm/133} {\path{doi:10.1090/gsm/133}}.
  \MR{2918544}

\bibitem[Ren08]{Rendall_book}
A.~Rendall, \emph{Partial {D}ifferential {E}quations in {G}eneral
  {R}elativity}, vol.~16, Oxford University Press, 2008.

\bibitem[Rob01]{robinson2001}
J.~Robinson, \emph{Infinite {D}imensional {D}ynamical {S}ystems}, Cambridge
  University Press, 2001.

\bibitem[Sog95]{Sogge_95}
C.~D. Sogge, \emph{Lectures on nonlinear wave equations}, Monographs in
  Analysis, II, International Press, Boston, MA, 1995. \MR{1715192}

\bibitem[Spe09]{Speck_09}
J.~Speck, \emph{Well-posedness for the {E}uler-{N}ordström system with
  cosmological constant}, J. Hyperbolic Differ. Equ. \textbf{6} (2009), no.~2,
  313--358, URL: \url{http://dx.doi.org/10.1142/S0219891609001885}. \MR{2543324
  (2011a:35529)}

\bibitem[SS98]{Shatah-Struwe-98}
J.~Shatah and M.~Struwe, \emph{Geometric wave equations}, Courant Lecture Notes
  in Mathematics, vol.~2, New York University, Courant Institute of
  Mathematical Sciences, New York; American Mathematical Society, Providence,
  RI, 1998. \MR{1674843}

\bibitem[Str86]{STRA2}
N.~Straumann, \emph{General {R}elativity and {R}elativistic {A}strophysics},
  Texts and Monographs in Physics, vol. 265, Springer Verlag, Heidelberg, 1986.

\bibitem[Tay97a]{taylor97}
M.~Taylor, \emph{Partial {D}ifferential {E}quations {I} {B}asic {T}heory},
  Springer, 1997.

\bibitem[Tay97b]{taylor97c}
\bysame, \emph{Partial {D}ifferential {E}quations {III} {N}onlinear
  {E}quations}, Springer, 1997.

\bibitem[TB01]{Todorova_Yordanov_01}
G.~Todorova and Y.~B., \emph{Critical exponent for a nonlinear wave equation
  with damping}, J. Differential Equations \textbf{174} (2001), no.~2,
  464--489.

\bibitem[Wan21]{Wang_2021}
X.~Wang, \emph{Propagation of regularity and long time behavior of 3{D} massive
  relativistic transport equation {I}: {V}lasov-{N}ordström system}, Comm.
  Math. Phys. \textbf{382} (2021), no.~3, 1843--1934, URL:
  \url{https://doi.org/10.1007/s00220-021-03987-2}, \href
  {http://dx.doi.org/10.1007/s00220-021-03987-2}
  {\path{doi:10.1007/s00220-021-03987-2}}. \MR{4232781}

\bibitem[Yag05]{Yagdjian_05}
K.~Yagdjian, \emph{Global existence in the {C}auchy problem for nonlinear wave
  equations with variable speed of propagation}, New Trends in the Theory of
  Hyperbolic Equations, Operator Theory: Advance and Applications, vol. 159,
  Birkhäuser, Basel, 2005, pp.~301--385, URL:
  \url{https://doi.org/10.1007/3-7643-7386-5_4}. \MR{2175919}

\end{thebibliography}
